\documentclass[journal]{IEEEtran}

\usepackage[utf8]{inputenc}
\usepackage{xcolor}
\usepackage{color, soul}
\usepackage{enumitem}
\setlist{nosep}

\usepackage[bottom]{footmisc}
\usepackage[hidelinks]{hyperref}
\usepackage{cite}
\usepackage{tabu}
\usepackage{arydshln}
\usepackage{diagbox}
\usepackage{bm, bbm}
\usepackage{mathtools}
\usepackage{nicefrac}
\usepackage{empheq}
\usepackage{mathtools}

\usepackage{accents}
\usepackage{units}
\usepackage{enumitem}

\usepackage{amsmath,amsfonts,amsthm,amssymb}
\usepackage[mathscr]{euscript}
\usepackage[bb=boondox]{mathalfa}
\usepackage[algoruled]{algorithm2e}
\usepackage{array}
\usepackage{mdframed}

\usepackage[noabbrev,capitalize]{cleveref}

\crefname{equation}{}{}
\Crefname{equation}{}{}
\crefname{section}{}{}
\Crefname{section}{}{}

\theoremstyle{definition} 

\newtheorem{condition}{Condition}

\theoremstyle{plain} 
\newtheorem{proposition}{Proposition}
\newtheorem{theorem}{Theorem} 
\newtheorem{definition}{Definition}

\theoremstyle{remark} 
\newtheorem{remark}{Remark}

\usepackage{tikz}

\newcommand{\set}[1]{\mathscr{#1}} 
\newcommand{\norm}[1]{\left\lVert#1\right\rVert} 
\newcommand{\dprod}[2]{\langle#1,#2\rangle}

\newcommand*{\tran}{^{\hspace{-0.1em}\mkern-1.mu\mathsf{T}}}

\DeclareMathOperator{\cvar}{CVaR}

\DeclareMathOperator{\Lap}{Lap}

\DeclareMathOperator{\dist}{dist}

\def\myproofstart{\noindent\textit{Proof. }}
\def\myproofend{\hspace*{\fill} $\square$ \vspace{0.3\baselineskip}}

\makeatletter
\newcommand{\pushright}[1]{\ifmeasuring@#1\else\omit\hfill$\displaystyle#1$\fi\ignorespaces}
\newcommand{\pushleft}[1]{\ifmeasuring@#1\else\omit$\displaystyle#1$\hfill\fi\ignorespaces}
\makeatother

\newcolumntype{C}[1]{>{\centering\arraybackslash}p{#1}} 

\newcounter{box}

\usepackage{placeins}
\newcommand{\subparagraph}{}
\usepackage{titlesec}

\title{Data Valuation from Data-Driven Optimization}

\author{Robert Mieth, Juan M.~Morales, H.~Vincent Poor}

\begin{document}

\pagestyle{empty}
\bstctlcite{IEEE:BSTcontrol} 

\maketitle

\begin{abstract}

With the ongoing investment in data collection and communication technology in power systems, data-driven optimization has been established as a powerful tool for system operators to handle stochastic system states caused by weather- and behavior-dependent resources.
However, most methods are ignorant to data quality, which may differ based on measurement and underlying privacy-protection mechanisms.
This paper addresses this shortcoming by (i) proposing a practical data quality metric based on Wasserstein distance, (ii) leveraging a novel modification of distributionally robust optimization using information from multiple datasets with heterogeneous quality to valuate data, (iii) applying the proposed optimization framework to an optimal power flow problem, and (iv) showing a direct method to valuate data from the optimal solution. 
We conduct numerical experiments to analyze and illustrate the proposed model and publish the implementation open-source.

\end{abstract}

\section{Introduction}

New means of data collection and processing can support power system operations by improving system observability and forecasting.  
For example, power system operators have widely adopted Phasor Measurement Units (PMUs) to better monitor the system state and are still investing in their adoption today \cite{pjm2020pmus}. 
On the distribution level, where traditional system planning has mostly relied on generously designed capacities and static protection schemes, thus creating a ``\textit{remarkably data-poor}''~\cite{siemens2019gridedge} system, emerging data collection is primarily driven by the ongoing smart meter roll-out. 
Additionally, third-party smart home devices such as \textit{Alexa} and \textit{Google Nest} collect high-fidelity grid-edge data on both energy consumption and consumers themselves. 
Yet, the data operationalization and exchange in the power and energy sector lags behind other industries, like finance and insurance, where a billion-dollar \textit{data economy}---characterized by comprehensive data collection, analysis, and monetization---has been widely adopted \cite{agarwal2019marketplace}.
On the one hand, this can be attributed to the high standards of critical power infrastructure on (cyber)security and customer privacy protection, which may be compromised by technology relying on the collection and exchange of data \cite{acharya2022false,giaconi2021smart}.
On the other hand, the complex nature of decision-making processes in a system governed by power system physics, complicate \textit{data valuation}, which, in turn, obstructs the design of (economic) frameworks and analyses for data collection, processing, and exchange. 
Motivated by the latter barrier, this paper proposes a data-to-decision approach leveraging distributionally robust optimization (DRO) that implicitly computes the value of available datasets in the context of their explanatory quality and the (physics-constrained) decision-making problem at hand.

Motivated by the need to deal with increasing uncertainty from weather-dependent (e.g., wind- and solar-powered generation) and behavior-dependent (e.g., EV charging and thermostatically controlled loads) energy resources, recent research highlights the usefulness of data-informed decision-making processes in power systems \cite{roald2023power}.
On the distribution (and electricity retail) level, for example, data streams collected from smart meters or other smart electric appliances have been shown to inform optimal control for smart loads \cite{hassan2020data,peng2022markovian}, leverage flexibility from distributed and controllable generators and batteries \cite{dall2017chance,mieth2018data,morales2021learning}, or enable optimal demand response pricing \cite{li2017distributed,mieth2019online,tucker2020constrained}.
On the transmission (and electricity wholesale) level, data-driven approaches have been shown to improve optimal generator dispatch \cite{guo2018data,arrigo2022wasserstein,morales2023prescribing}, system resiliency \cite{crozier2022data}, and energy forecasting and trading \cite{fernandez2021inverse,munoz2020feature}. 

The approaches proposed in \cite{hassan2020data,peng2022markovian,dall2017chance,mieth2018data,li2017distributed,mieth2019online,morales2021learning,morales2023prescribing,tucker2020constrained,guo2018data,arrigo2022wasserstein,crozier2022data,fernandez2021inverse,munoz2020feature} assume that the required data streams are freely available.
However, collecting and sharing data can  create costs either in the form of direct costs related to data storage, processing, and transmission \cite{ren2018datum} or indirect costs related to a loss of privacy \cite{bessa2018data} or competitive disadvantages \cite{goncalves2020towards}. 
Moreover, the monetary value of data access is fully absorbed by the data user \cite{han2021monetizing} and the resulting inequitable allocation of costs and benefits of data exchange creates smart grid acceptance barriers \cite{le2020ethical}.

To overcome this disparity and enable transparent data valuation, recent research in the power and energy community has turned its attention to markets for data \cite{goncalves2020towards,han2022trading,pinson2022regression,han2021monetizing} and forecasts \cite{raja2023market,xie2022robust}. 
Proposals for data markets in energy systems, as in \cite{goncalves2020towards,pinson2022regression}, largely follow the design proposed in \cite{agarwal2019marketplace}, where a central trustworthy entity collects prediction tasks from data users and matches them with datasets and cost bids from data owners.
While this approach has been shown to incentivize data owners to share data by creating a new revenue stream, their application has so far been limited to regression-based forecasting markets \cite{goncalves2020towards,han2022trading,pinson2022regression}. 
Here, and similarly in \cite{han2021monetizing}, datasets are valuated based on their Shapely Value, i.e., by computing the performance of the prediction task with and without each dataset and evaluating the difference.
While this design has the advantage of only sharing data with the data market and not with the individual data users, data valuation hinges on the computationally expensive Shapely value computation (or its approximation), is sensitive to how datasets are split into testing and training data, and requires a new centralized service. 

In this paper we develop an alternative and more direct data valuation approach that is in line with proposals like \cite{hassan2020data,peng2022markovian,dall2017chance,mieth2018data,li2017distributed,mieth2019online,morales2023prescribing,guo2018data,arrigo2022wasserstein,crozier2022data,morales2021learning,tucker2020constrained,fernandez2021inverse,munoz2020feature}, where a decision
maker (e.g., a power system operator or utility) has access to datasets from multiple sources to inform one or multiple uncertain parameters. 
These datasets may have been subject to privacy protection mechanisms, e.g., as discussed in \cite{giaconi2021smart,vdvorkin2020differentially}, or noisy measurements, e.g., as is common in PMUs \cite{brown2016characterizing}, which may reduce the data quality, i.e., its usefulness for the intended decision-making task.
Considering this context, we propose a variant of data-driven DRO that implicitly computes marginal data value as a function of data quality.
We introduce this method by making the following contributions:
\begin{enumerate}[leftmargin=1.4em, labelwidth=1.4em]
\item We {\color{black} propose a metric that quantifies the quality of a dataset in terms of its ability to infer the distribution of an unknown parameter. We argue for the suitability of the Wasserstein distance as a basis for this metric and discuss connections to common data obfuscation methods.} 
\item {\color{black} We introduce the idea of \textit{marginal data value}. For this, we} derive a novel modification of data-driven DRO with Wasserstein ambiguity sets, a popular technique in general operations research \cite{esfahani2018data} and smart grid applications \cite{guo2018data}. 
In line with recent studies on DRO with fragmented datasets \cite{awasthi2022distributionally}, our DRO modification allows the use of multiple datasets with heterogeneous data quality. 
To our knowledge, this is the first time that the rationale behind DRO is exploited, not for \emph{data-driven decision making}, but for \emph{decision-driven data valuation}. 
{\color{black} Importantly, our approach can be applied to value data for decision-making problems that are not necessarily intended to be solved by DRO.}
\item We apply the proposed method to an optimal power flow (OPF) problem with data-informed uncertain parameters and discuss how data value can be directly derived from the problem. 
We present a numerical case study to illustrate and discuss the proposed method.
\end{enumerate}
By deriving marginal data value directly from the decision-making problem and its constraints (e.g., power flow physics), we obtain an interpretable quantity that aligns with established marginal cost-based valuation approaches in power economics.
{\color{black} This approach departs from existing data valuation concepts that rely on a market mechanism \cite{goncalves2020towards,han2022trading,pinson2022regression,han2021monetizing,raja2023market,xie2022robust} by providing a marginal data value that is implicit to the decision-making process and does not require a secondary computation (e.g., to compute Shapely value) or a preliminary bidding process.
}

\section{Quantifying Data Quality}
\label{sec:privacy_and_dist_ambiguity}

As a first step, {\color{black} we take the perspective of data owners and establish a metric to quantify the quality of a dataset.
We consider $D$ data owners indexed by $j=1,...,D$ and assume that each data owner $j$ has access to a dataset $\set{X}_j$ that contains information on an uncertain variable (\textit{feature}) $\bm{\xi}_j$.
}
For ease of notation, but without loss of generality, we assume a one-to-one relationship between data owners, datasets, and features so that all can be uniquely identified by index $j$.
{\color{black}
We consider each feature $\bm{\xi}_j$ to be a quantity of interest for a decision maker, e.g., future demand, renewable power injections, or state estimations. 
The data in dataset $\set{X}_j$ can be used to estimate $\bm{\xi}_j$ and, because the true value of $\bm{\xi}_j$ may only be observed in hindsight or not at all, to make inferences on its distribution $\mathcal{P}_j$.
For this purpose,
}
dataset $\set{X}_j$ may, for example, contain smart meter measurements, power generation measurements with auxiliary weather data from a commercial wind or solar farm, or PMU data.
{\color{black}
For the purpose of this paper, we therefore establish the following definition:
\begin{definition}[Data quality]\label{def:data_quality}
   The quality of a dataset $\set{X}$ with respect to a random variable $\bm{\xi}$ is a measure of how accurately the distribution $\mathcal{P}$ of $\bm{\xi}$ can be estimated using the data from~$\set{X}$.
\end{definition}
Section~\ref{ssec:data_quality_measure} below formalizes this definition using the Wasserstein distance and provides necessary discussion. 
The following Section~\ref{ssec:modeling_obfuscation_and_noise} then describes how to model noise and intentional data obfuscation (e.g., for privacy protection) in this framework.

\subsection{Data quality as distributional ambiguity}
\label{ssec:data_quality_measure}

Consider $\mathcal{P}_j$ as the true distribution of $\bm{\xi}_j$ and consider $\set{X}_j$ as all data that is available to compute an estimate $\widehat{\mathcal{P}}_j$ of $\mathcal{P}_j$. 
If the information contained in the data of $\set{X}_j$ is sufficient and allows an estimate $\widehat{\mathcal{P}}_j$ close to $\mathcal{P}_j$, i.e., 
a distance $\dist\!\big(\widehat{\mathcal{P}}_j, \mathcal{P}_j\big)$ is small, we say that dataset $\set{X}_j$ is of high quality. Conversely, if the best possible estimate $\widehat{\mathcal{P}}_j$ that $\set{X}_j$ allows is likely to be very different from $\mathcal{P}_j$, i.e., $\dist\!\big(\widehat{\mathcal{P}}_j, \mathcal{P}_j\big)$ is potentially large, we say that dataset $\set{X}_j$ is of low quality. 
}

A natural and practical candidate for
{\color{black}
$\dist\!\big(\widehat{\mathcal{P}}_j, \mathcal{P}_j\big)$ 
is the \textit{Wasserstein distance}. For a given pair $\widehat{\mathcal{P}}_j$ and $\mathcal{P}_j$,
}
their \mbox{$p$-Wasserstein} distance to the power of $p$ is defined as:
\begin{equation}
    \mathcal{W}_p^p(\mathcal{P}_j, \widehat{\mathcal{P}}_j) = \inf_{\mathcal{J} \in \set{J}(\mathcal{P}_j, \widehat{\mathcal{P}}_j)} \int \norm{\bm{\xi}-\widehat{\bm{\xi}}}^{\color{black}p} d \mathcal{J}(\bm{\xi},\widehat{\bm{\xi}}),
\label{eq:wasserstein_definition}
\end{equation}
where $\set{J}(\mathcal{P}_j, \widehat{\mathcal{P}}_j)$ is the set of all joint distributions $\mathcal{J}$ with marginals $\mathcal{P}_j$ and $\widehat{\mathcal{P}}_j$, and $\norm{\cdot}$ is  a norm. 
{\color{black}
Intuitively, the Wasserstein distance}
computes the minimal amount of probability mass that must be ``moved'' to construct $\widehat{\mathcal{P}}_j$ from~$\mathcal{P}_j$. 

{\color{black}
The Wasserstein distance possesses desirable properties that motivate its selection for our data quality metric in the context of data-driven decision-making processes.
In particular, if $\widehat{\mathcal{P}}_j$ is an empirical distribution supported by the data in $\set{X}_j$, then}
the Wasserstein distance is particularly descriptive when the difference between $\widehat{\mathcal{P}}_j$ and $\mathcal{P}_j$ is
{\color{black}
driven by random noise that perturbs the data in $\set{X}_j$ \cite{esfahani2018data}. This is a common scenario in practical settings, which we discuss in  Section~\ref{ssec:modeling_obfuscation_and_noise} below.
Indeed, in this case other popular dissimilarity measures are inexpressive. 
For example, the Kullback-Leibler divergence is unable to compare distributions with disjoint support, which is true almost surely for the empirical distribution supported on $\set{X}_j$ and its noisy counterpart.
In this case, the Kullback-Leibler divergence is defined as infinite.
Similarly, the total variation distance will equal one almost surely \cite[p.~610]{alvarez2012similarity} in this setting.
In contrast, the $p$-Wasserstein metric is well-behaved when working with empirical distributions and exhibits asymptotic consistency properties \cite[Th.~1]{alvarez2012similarity}.
Moreover, it allows a tractable integration in decision-making problems, as we discuss in Section~\ref{sec:multi-source_dro} below.
These properties have contributed the the popularity of the Wasserstein metric for data-driven optimization \cite{esfahani2018data} and machine learning \cite{kuhn2019wasserstein}, and also motivates our choice here.

In many practical settings the true distribution of $\bm{\xi}_j$ is unknown and $\mathcal{W}_p^p(\mathcal{P}_j, \widehat{\mathcal{P}}_j)$ is a random variable that depends on the finite samples available in $\set{X}_j$ and any possible perturbation these samples may have suffered due to contamination, noise or intentional obfuscation.
We therefore use a \textit{confidence upper bound} $\epsilon_j$ on $\mathcal{W}_p^p(\mathcal{P}_j, \widehat{\mathcal{P}}_j)$ as a suitable proxy for data quality. This upper bound can be readily estimated by the data owner \cite{esfahani2018data,duan2018distributionally} and allows for useful extensions in practical settings with perturbed data, which we discuss in the following section.

\subsection{Practical considerations on data quality and obfuscation}
\label{ssec:modeling_obfuscation_and_noise}

In practical applications, dataset $\set{X}_j$ may neither be equal to the raw dataset that data owner $j$ is collecting from one or multiple data sources, nor may the data owner wish to publish $\set{X}_j$ directly to protect sensible information. 
In the following we describe a model that relates $\set{X}_j$ to ``raw'' and ``publishable'' datasets, discuss the role of data quality using the concept from Section~\ref{ssec:data_quality_measure} above, and provide three practical examples.

We denote the raw dataset, i.e., all data that is available to data owner $j$ as $\set{Y}_j$. This data may contain sensitive information that must be protected.
Dataset $\set{Y}_j$ is transformed into dataset $\set{X}_j$ through a \textit{query} $\mathcal{Q}_j$.}
Query $\mathcal{Q}_j$ captures all data manipulations required for making the raw dataset $\set{Y}_j$ useful for the intended application, i.e., inferences on feature $\bm{\xi}_j$. This may include aggregation, histogram computation, forecasting, or other numerical outputs that have full access to $\set{Y}_j$ to create $\set{X}_j$ and make it \textit{actionable} for the intended purpose and ensure the \textit{interoperability} between data collector and data recipient \cite{gopstein2021nist}. 
{\color{black}%
Per definition, query $\mathcal{Q}_j$ conserves all information on feature $\bm{\xi}_j$ available in $\set{Y}_j$.}
However, it has been shown that many queries, including aggregation and forecasting, can be reverse-engineered to reveal specific information on the data source, e.g., presence of certain appliances, home occupancy, or production schedules \cite{giaconi2021smart,molina2010private,hassan2021privacy}.

\begin{figure}
    \centering
    \includegraphics[width=0.95\linewidth]{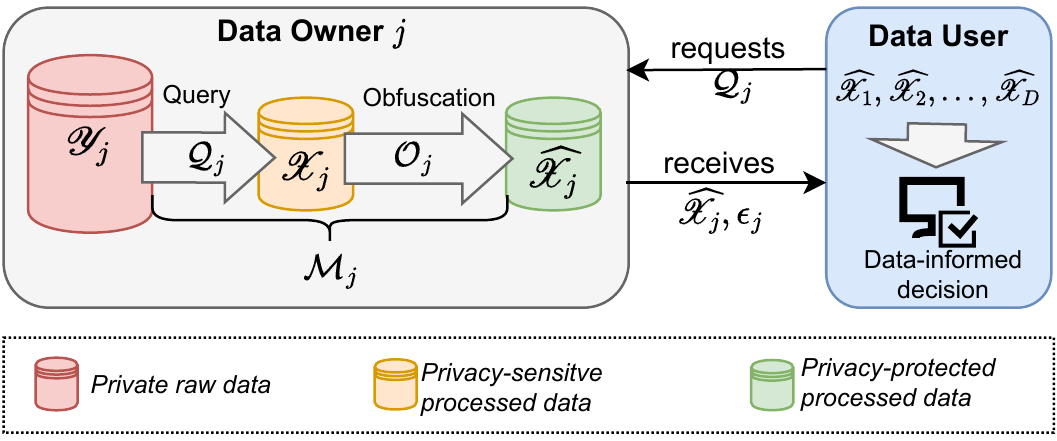}
    \caption{Data owner $j$ processes ($\mathcal{Q}_j$) and obfuscates ($\mathcal{O}_j$) data based on the query of the data user. Data user receives datasets $\{\widehat{\set{X}}_j\}_{j=1}^D$ alongside their quality signal $\{\epsilon_j\}_{j=1}^D$ for its intended data-informed decision-making.}
    \label{fig:mechanism_scheme}
\end{figure}

{\color{black}
To protect sensitive information contained in $\set{X}_j$, the data owners may apply an \textit{obfuscation} $\mathcal{O}_j$ to  $\set{X}_j$.}
Obfuscation $\mathcal{O}_j$ alters $\set{X}_j$ into a new dataset $\widehat{\set{X}}_j$. This dataset retains the properties of $\set{X}_j$ that are required by the data user (e.g., resolution or data format), but additionally meets privacy requirements imposed by the data owner. 
Common obfuscation approaches are differential privacy \cite{dwork2014algorithmic} or algorithmic obfuscation techniques \cite{giaconi2021smart}. 
{\color{black}%
We call the composition of $\mathcal{Q}_j$ and $\mathcal{O}_j$ a \textit{mechanism} $\mathcal{M}_j$, leading to the following formal relationship between the datasets:
\begin{align}
    \widehat{\set{X}}_j = \mathcal{O}_j(\set{X}_j) = (\mathcal{O}_j \circ \mathcal{Q}_j)(\set{Y}_j) = \mathcal{M}_j(\set{Y}_j).
\end{align}
Fig.~\ref{fig:mechanism_scheme} illustrates this relationship.
Notably, applying obfuscation $\mathcal{O}_j$ impacts the data quality of $\widehat{\set{X}}_j$ relative to $\set{X}_j$ as per Definition~\ref{def:data_quality}.
The introduction of obfuscation $\mathcal{O}_j$ will reduce the information on $\bm{\xi}_j$ contained in $\widehat{\set{X}}_j$ such that the data quality bound $\epsilon_j$ increases.
The following examples discuss common obfuscations and their impact on data quality.}

\subsubsection{Example: Differential privacy}
\label{ssec:differential_privacy}

In differential privacy, obfuscation $\mathcal{O}_j$ is an additive noise drawn from a parametrized distribution such that $\mathcal{M}^{\rm DP}$ is
\begin{equation}
    \mathcal{M}^{\rm DP} = \mathcal{Q}_j(\set{Y}_j) + \mathcal{O}_j(\theta, \mathcal{Q}_j, \set{Y}_j),
\end{equation}
where the additonal parameter $\theta$ determines the level of desired privacy.
For example, to achieve $\theta$-differential privacy using the popular Laplace Mechanism \cite{dwork2014algorithmic,vdvorkin2020differentially}, each element of $\set{X}_j = \mathcal{Q}_j(\set{Y}_j)$ is {perturbed} by an i.i.d. sample from a Laplace distribution $\Lap(\Delta \mathcal{Q}_j/\theta)$, where $\Delta \mathcal{Q}_j$ defines the worst-case sensitivity of $\mathcal{Q}_j$ to changing exactly one element in $\set{Y}_j$ \cite{dwork2014algorithmic}.
Denoting $X_j$ and $\widehat{X}_j$ as random variables described by distributions $\mathcal{P}_j$ and $\widehat{\mathcal{P}}_j$, respectively, and $Z_j$ as the random variable described by distribution $\Lap(\Delta \mathcal{Q}_j/\theta)$, we can express the information loss caused by obfuscation in $\widehat{\set{X}}_j$ relative to $\set{X}_j$ in terms of Wasserstein distance
\begin{align}
    \mathcal{W}_p^p(X_j, \widehat{X}_j) 
        &= W_p^p(X_j, X_j + Z_j) \nonumber \\
        &\stackrel{\text{(A)}}{\le} \mathcal{W}_p^p(X_j, X_j) + \mathcal{W}_p^p(0, Z_j) \stackrel{\text{(B)}}{=} \mathcal{W}_p^p(0, Z_j) \nonumber \\
        & \stackrel{\text{(C)}}{=} \mathbb{E}\norm{Z_j}^p,
\label{eq:wasserstein_additive}
\end{align}
where relation (A) uses the additive property of the Wasserstein distance \cite{panaretos2019statistical}, relation (B) uses the fact that the Wasserstein distance between two identical distributions is zero, and relation (C) results from the definition of the Wasserstein distance.
Note that in \cref{eq:wasserstein_additive} we extend the $\mathcal{W}_p$ notation introduced in \cref{eq:wasserstein_definition} by using $X_j\sim\mathcal{P}_j$ and $\widehat{X}_j\sim\widehat{\mathcal{P}}_j$ as arguments. 
The derivation in \cref{eq:wasserstein_additive} shows that for additive noise on the original data the $p$-Wasserstein distance is bounded from above by $\mathbb{E}\norm{Z_j}^p$, which allows for a computationally efficient direct approximation of the upper bound of the Wasserstein distance. 
For $p=1$ and and an underlying 1-norm this is equal to the expectation of the absolute value of $Z_j$. For $p=2$ and an underlying 2-norm this is equal to the variance of $Z_j$.

\begin{remark}
\label{rem:publishing_epsilon}
Publishing $\mathbb{E}\norm{Z_j}^p$ alongside $\widehat{\set{X}}_j$ does not compromise differential privacy as both the query $\mathcal{Q}_j$ \textit{and} obfuscation $Z_j$ are considered public information \cite{wasserman2010statistical}.
\end{remark}

\subsubsection{Example: Noisy data}
\label{ssec:noisy_data}
The dataset $\set{Y}_j$ may not contain privacy-relevant data but can be subject to noisy measurements, e.g., time-series from utility-owned PMUs.
In this case, we can use obfuscation $\mathcal{O}_j$ to model this noise. We set $\mathcal{Q}_j(\set{Y}_j)=\set{Y}_j$ such that $\set{X}_j\equiv\set{Y}_j$ and $\mathcal{O}_j$ such that it adds noise to each data entry in $\set{X}_j$. If the noise can be assumed i.i.d for each sample, as suggested by empirical research on PMUs in \cite{brown2016characterizing}, we recover the result from \cref{eq:wasserstein_additive}. 
Hence, \cref{eq:wasserstein_additive} allows to bound the Wasserstein distance even if the  original data is unknown and only information on the noise is available.

\subsubsection{Example: Aggregation protocols}
\label{ssec:algorithmic_obfuscation}

Instead of obfuscating the entries of $\set{X}_j$ with random noise as in differential privacy, $\mathcal{O}_j$ may follow a more elaborate protocol. For example, aggregation protocols for smart meters add masking values to individual smart meter measurements that cancel each other out either over time (e.g., to ensure that the total energy consumption for a billing period is exact) or over multiple data sources (e.g., to ensure that the total energy consumption for a group of households is exact) \cite{giaconi2021smart,kursawe2011privacy}. In this case, the approximation in \cref{eq:wasserstein_additive} may no longer hold and the data owner has to compute a bound on the Wasserstein distance between the distribution supported by $\set{X}_j$ and $\widehat{\set{X}}_j$ as per \cref{eq:wasserstein_definition}, which is a tractable computation.

\section{Multi-Source data-driven DRO}
\label{sec:multi-source_dro}

In this section we take the perspective of the data user (see Fig.~\ref{fig:mechanism_scheme}), e.g., a power system operator, who receives datasets $\{\widehat{\set{X}}_j\}_{j=1}^D$ alongside their quality information, {\color{black} i.e., confidence upper bounds $\{\epsilon_j\}_{j=1}^D$} of the Wasserstein distance to the true distribution as discussed in Section~\ref{sec:privacy_and_dist_ambiguity} above. {\color{black} In the following, we will refer to  $\epsilon_j$ as data quality for conciseness.}
We now assume that the data user wishes to to make a decision $\bm{x}$ such that a cost function $c(\bm{x}, \bm{\xi}_1,...,\bm{\xi}_D)$ that depends on decision $\bm{x}$ and on the $D$ uncertain parameters $\bm{\xi}_1,...,\bm{\xi}_D$ is minimized.   
{\color{black}%
The resulting problem takes the general form 
\begin{equation}
    \inf_{\bm{x}\in\set{F}} \mathbb{E}_{\mathcal{Q}}[c(\bm{x}, \bm{\xi}_1,...,\bm{\xi}_D)],
\label{eq:general_dec}
\end{equation}
where $\set{F}$ is a feasible region constraining the vector of decision variables, $\mathcal{Q}$ is a joint distribution of the random vectors $\bm{\xi}_1,...,\bm{\xi}_D$, and $\mathbb{E}_{\mathcal{Q}}$ may be the expectation operator or any suitable risk measure. We focus on expectation for the remainder of the paper and point to \cite{mieth2020risktrading} for further discussion.
}

The data user can extract information on the distribution of each $\bm{\xi}_j$ from $\widehat{\set{X}}_j$ \textit{with respect to data quality} $\epsilon_j$. 
{\color{black}%
A natural reformulation of \cref{eq:general_dec} that internalizes such information 
} is a data-driven DRO problem.
{\color{black}%
Data-driven DRO has gained popularity in many domains including power system operations and planning (e.g., \cite{mieth2018data,arrigo2022embedding,hassan2019stochastic,guo2018data,esteban2021distributionally,roald2023power}) as it provides an optimal balance between robustness and optimism in data-driven optimization \cite{van2021data}.
However, none of the existing DRO formulations internalize distributional ambiguity information \textit{individually} for each feature $j$ given $\widehat{\set{X}}_j$ and $\epsilon_j$.
In the remainder of this section, we bridge this gap and introduce two practical conditions that improve the tractability of the general formulation.}

The general DRO reformulation of \cref{eq:general_dec} is
\begin{equation}
    \inf_{\bm{x}\in\set{F}} \sup_{\mathcal{Q}\in\set{A}} \mathbb{E}_{\mathcal{Q}}[c(\bm{x}, \bm{\xi}_1,...,\bm{\xi}_D)],
\label{eq:general_dro}
\end{equation}
{\color{black}%
where, in contrast to \cref{eq:general_dec}, distribution} 
$\mathcal{Q}$ is drawn from a set of candidate distributions $\set{A}$ (called ambiguity set), {\color{black} which is informed by $\{\widehat{\set{X}}_j\}_{j=1}^D$ and $\{\epsilon_j\}_{j=1}^D$.}
From each dataset $\widehat{\set{X}}_j$ the data user can construct the empirical distribution $\widehat{\mathcal{P}}_j$. With the additional information that the distance between $\widehat{\mathcal{P}}_j$ and the true distribution $\mathcal{P}_j$ (see also Section~\ref{sec:privacy_and_dist_ambiguity}) is bounded by $\epsilon_j$ with high confidence, we can define the \textit{multi-source Wasserstein ambiguity set}
\begin{equation}
    \set{A}^{\rm MSW} \!\!\coloneqq\! 
   \left\{ \mathcal{Q}\in\set{P}(\Xi) \left|
    \begin{array}{ll}
          P_{j\#}\mathcal{Q} = \mathcal{Q}_j, & j=1,...,D  \\
         \mathcal{W}_p^p(\widehat{\mathcal{P}}_j, \mathcal{Q}_j) \le \epsilon_j, &  j=1,...,D 
    \end{array}\right.\right\}, 
\label{eq:multi_source_ambiguity_set}
\end{equation}
where $\set{P}(\Xi)$ is the set of all possible probability measures with support $\Xi$. Note that the support $\Xi$ can be inferred from the available datasets or based on some other technical considerations (see discussion in Section~\ref{ssec:data_and_support} below). Further, $P_{j\#}$ denotes the push-forward distribution of the joint measure $\mathcal{Q}$ under the projection onto the $j$-th coordinate, i.e.,
\begin{equation}
    P_{j\#}\mathcal{Q}(d\bm{\xi}_j^*)\!\coloneqq\!\! \int_{\Xi_{-j}} \hspace{-0.3cm}\mathcal{Q}(d\bm{\xi}_1,...,d\bm{\xi}_{j-1}, d\bm{\xi}_j^*, d\bm{\xi}_{j+1}, ..., d\bm{\xi}_D)
\end{equation}
with $\Xi_{-j} \coloneqq \Xi_1 \times ... \times \Xi_{j-1} \times \Xi_{j+1} ... \times \Xi_{D}$ and $\Xi_j$ being the projection of $\Xi$ on the $j$-th coordinate.
The ambiguity set $\set{A}^{\rm MSW}$ generalizes established Wasserstein ambiguity sets, e.g., as in \cite{esfahani2018data}, to accommodate individual Wasserstein budgets $\epsilon_j$ for each feature $j$.

\begin{theorem}
\label{th:multi-source-wasserstein-dro}
A DRO problem as in \cref{eq:general_dro} with ambiguity set $\set{A}=\set{A}^{\rm MSW}$ as in \cref{eq:multi_source_ambiguity_set} can be rewritten as
\allowdisplaybreaks
\begin{subequations}
\begin{align}
\inf_{\substack{\bm{x}\in\set{F},\\ \lambda_j\ge0,s_{\iota}}}\ & \sum_{j=1}^D \lambda_j \epsilon_j + \frac{1}{|\set{I}|}\sum_{\iota \in \set{I}} s_{\iota} \\
\text{s.t.}\quad & s_{\iota} \ge \sup_{\bm{\xi}\in\Xi}c(\bm{x},\bm{\xi}) - \sum_{j=1}^D\lambda_j\!\norm{\bm{\xi}_j - \widehat{\bm{\xi}}_{\iota_j}\!}^p,\  \forall \iota \in \set{I}
\end{align}%
\label{eq:mwsdro_general_formulation}%
\end{subequations}%
\allowdisplaybreaks[0]%
where $\iota \in \set{I} = \bigtimes_{j=1}^D\{1,\ldots,N_j\}$ is a multi-index from the set of all combinations of data-sample indices across all datasets, $N_j$ denotes the number of data points in set $\widehat{\set{X}}_j$, i.e., $|\widehat{\set{X}}_j|$, and $\lambda_j$, $s_{\iota}$ are auxiliary variables.
Note that we use $\bm{\xi}$ (without index) as shorthand for dependency on $\bm{\xi}_1,...,\bm{\xi}_D$.
\end{theorem}
\myproofstart
See Appendix~\ref{ax:proof_of_th_wassertein_dro_general}.
\myproofend

The formulation in \cref{eq:mwsdro_general_formulation} makes no assumptions on how the samples of the different datasets may be related, thus capturing all possible dependencies between the empirical distributions of the individual datasets. As a result, the number of constraints grows exponentially with the number of features $D$ (i.e., the number of constraints is $\prod_{j=1}^D N_j$).
However, under \textit{either} one of two conditions with practical relevance, the problem complexity can be dramatically reduced. 

\renewcommand*{\thecondition}{C-Sep.}
\begin{condition}[Separable cost function] 
\label{cond:separable_cost_and_support}
The cost function $c(\bm{\xi}_j)$ can be written as $\sum_{j=1}^Dc_j(\xi_j)$, i.e., it is separable in the random vectors associated with each data provider.
\end{condition}

\renewcommand*{\thecondition}{C-Std.}
\begin{condition}[Standardized data]
\label{cond:standardized_data}
The datasets are standardized, i.e., they have the same lengths ($N_j=N',\ j=1,...,D$) and we can define a shared index $i'$ such that all samples $\{\widehat{\bm{\xi}}_{1i'},...,\widehat{\bm{\xi}}_{Di'}\}_{i'=1}^{N'}$ can be related, e.g., in terms of a shared time-stamp.
\end{condition}

\begin{proposition}\label{prop:mwsdro_sep_cost}
Under Condition~\ref{cond:separable_cost_and_support} the multi-source Wasserstein DRO from \cref{eq:mwsdro_general_formulation} can be reformulated as 
\begin{subequations}
\begin{align}
    \inf_{\substack{\bm{x}\in\set{F},\\ \lambda_j\ge0,s_{ji}}}\ & \sum_{j=1}^D \lambda_j \epsilon_j + \sum_{j=1}^D \frac{1}{N_j}\sum_{i=1}^{N_j} s_{ji} \label{prop1:objective} \\
    \text{s.t.}\quad & s_{ji} \ge \sup_{\bm{\xi}_j\in\Xi_j}\ c_j(\bm{x},\bm{\xi}_j) -  \lambda_j\!\norm{\bm{\xi}_j - \widehat{\bm{\xi}}_{ji}}_p, \\
    & \hspace{3.5cm} j=1,...,D,\ i=1,...,N_j \label{prop1:constraint}\nonumber
\end{align}%
\label{eq:mwsdro_general_formulation_sep_cost}%
\end{subequations}%
which scales linearly in $D$ and $N_j$.
\end{proposition}
\myproofstart
    See Appendix~\ref{ax:proof_sep_cost}.
\myproofend

\begin{proposition}
\label{prop:mwsdro_stand_data}
Under Condition~\ref{cond:standardized_data} the multi-source Wasserstein DRO from \cref{eq:mwsdro_general_formulation} can be reformulated as 
\allowdisplaybreaks
\begin{subequations}
\begin{align}
\inf_{\lambda_j\ge0}\ 
    & \sum_{j=1}^D \lambda_j \epsilon_j + \frac{1}{N'}\sum_{i'=1}^{N'} s_{i'} \\
\text{s.t.} \quad 
    & s_{i'} \ge \sup_{\bm{\xi}\in\Xi}\ c(\bm{x}, \bm{\xi}) - \sum_{j=1}^D\lambda_j\norm{\bm{\xi}_j - \widehat{\bm{\xi}}_{ji'}}_p,\ i'=1,...,N'
\end{align}%
\label{eq:mwsdro_general_formulation_stand_data}
\end{subequations}%
\allowdisplaybreaks[0]%
which scales linearly in $D$ and $N'$.
\end{proposition}
\myproofstart
See Appendix~\ref{ax:proof_stand_data}.
\myproofend

The result in Proposition~\ref{prop:mwsdro_stand_data} has a similar structure to established Wasserstein DRO, e.g.,\cite[Theorem 4.2]{esfahani2018data} where a single Wasserstein distance $\epsilon$ is used. 
We note that if $\sum_{j=1}^D\epsilon_j=\epsilon$, Proposition~\ref{prop:mwsdro_stand_data} will provide a less conservative worst-case expectation, because it considers ambiguity information for each marginal distribution in $j$ individually.
{\color{black}%
Finally, we highlight that the method presented in this section is applicable to any optimization problem with expectation constraints as in \cref{eq:general_dec} that uses data from multiple sources to inform uncertain features $\bm{\xi}_j,\ j=1,...,D$. 
}

\section{Multi-source data-driven OPF}
\label{sec:ms_opf}

We now apply the multi-source data-driven DRO formulation from Section~\ref{sec:multi-source_dro} above to an OPF problem with stochastic (e.g., weather- or behavior-dependent) resources.
We model the $D$-dimensional vector of power injections from these resources as $\bm{u}(\bm{\xi}) = \bm{u} + \bm{\xi}$, where $\bm{u}$ is a (deterministic) forecast and $\bm{\xi}$ is a vector of random forecast errors. 
Real-time power imbalance caused by errors $\bm{\xi}$ is balanced by $G$ controllable resources with output $\bm{p}(\bm{\xi}) = \bm{p} - \bm{A}\bm{\xi}$ such that $\bm{A}\tran\bm{1}_{G} = \bm{1}_D$, where $\bm{1}_{X}$ is a vector of ones of dimension $X$. This ensures that for each uncertain resource the total error is balanced \cite{mieth2019risk}.

\subsection{General formulation}
\label{ssec:general_model_formulation}

The system operator aims to decide on the power injection $\bm{p}$ alongside up- and downward reserve capacity ($\bm{r^+}$, $\bm{r^-}$) such that system cost are minimized and security constraints are maintained with high probability. 
We formulate this goal as a distributionally robust OPF problem with joint chance-constraints that protects the decision against uncertain injections $\bm{u}(\bm{\xi})$ and $\bm{p}(\bm{\xi})$:
\allowdisplaybreaks
\begin{subequations}
\begin{align}
\min_{\substack{\bm{p}, \bm{A}, \bm{r}^+, \bm{r}^-,\\ \bm{f}^{\rm RAM+}, \bm{f}^{\rm RAM-}}}  
    \hspace{-0.5cm}&\hspace{0.5cm} \dprod{\bm{c}^{\rm E}\!}{\bm{p}} \!+ \!\dprod{\bm{c}^{\rm R}\!}{\bm{r}^+\! + \bm{r}^{-}} \!+\!\! \max_{\mathcal{Q}\in\set{A}} \!\!\mathbb{E}_{\mathcal{Q}}[\dprod{\bm{c}^{\rm A}\!}{-\bm{A}\bm{\xi}}] \nonumber \\
\text{s.t.}\quad \label{base_dcopf_dro:objective} \\
(\pi):\quad    
    & \dprod{\bm{1}_G}{\bm{p}}  = \dprod{\bm{1}_V} {\bm{d}} - \dprod{\bm{1}_D}{\bm{u}} \label{base_dcopf_dro:enerbal}\\
(\bm{\chi}):\quad    
    & \bm{A}\tran\bm{1}_G= \bm{1}_D \label{base_dcopf_dro:resbal}\\
(\bm{\sigma^{\rm up}}):\quad     
    & \bm{p} + \bm{r}^+ \le \bm{p}^{\rm max} \label{base_dcopf_dro:gen_uplim}\\
 (\bm{\sigma^{\rm lo}}):\quad    
    & \bm{p} - \bm{r}^- \ge \bm{p}^{\rm min} \label{base_dcopf_dro:gen_dnlim}\\
(\bm{\beta^{\rm up}}):\quad   
    & \bm{B}^{\rm G} \bm{p} + \bm{B}^{\rm W} \bm{u} - \bm{B}^{\rm B}\bm{d} = \bm{f}^{\rm max} - \bm{f}^{\rm RAM+} \label{base_dcopf_dro:lin_uplim}\\
(\bm{\beta^{\rm lo}}):\quad      
    &-(\bm{B}^{\rm G} \bm{p} + \bm{B}^{\rm W} \bm{u} - \bm{B}^{\rm B}\bm{d}) = \bm{f}^{\rm max} - \bm{f}^{\rm RAM-} \label{base_dcopf_dro:lin_dnlim} \\
    &\inf_{\mathcal{Q}\in\set{A}} \mathcal{Q}\!\begin{Bmatrix*}[c]
        -\bm{A}\bm{\xi} \le \bm{r}^+ \\
        \bm{A}\bm{\xi} \le \bm{r}^-  \\
        (\bm{B}^{\rm W} - \bm{B}^{\rm G}\bm{A})\bm{\xi} \le \bm{f}^{\rm RAM+}\!  \\
        -(\bm{B}^{\rm W} - \bm{B}^{\rm G}\bm{A})\bm{\xi} \le \bm{f}^{\rm RAM-}\! 
    \end{Bmatrix*}\! \ge\! 1\!-\!\gamma \label{base_dcopf_dro:joint_cc}\\
    &\bm{p}, \bm{A}, \bm{r}^+, \bm{r}^-, \bm{f}^{\rm RAM+}, \bm{f}^{\rm RAM-} \ge 0, \label{base_dcopf_dro:nonneg}
\end{align}%
\label{eq:base_dcopf_dro}%
\end{subequations}%
\allowdisplaybreaks[0]%
where we use $\dprod{\cdot}{\cdot}$ to denote the dot-product of two vectors. Greek letters in parentheses denote dual multipliers.

The objective~\cref{base_dcopf_dro:objective} minimizes system cost defined by vectors for energy cost, reserve provision cost, and reserve activation cost, $\bm{c}^{\rm E}$, $\bm{c}^{\rm R}$, and $\bm{c}^{\rm A}$, respectively. Uncertain reserve activation $-\bm{A}\bm{\xi}$ is accounted for in terms of its worst-case expectation with respect to all distributions in the ambiguity set $\set{A}$.
Energy balance \cref{base_dcopf_dro:enerbal} ensures that the total generator injections equals the total system demand at all system buses minus the forecast injections from uncertain resources.
Similarly, \cref{base_dcopf_dro:resbal} ensures that  all forecast errors are balanced. 
Constraints \cref{base_dcopf_dro:gen_uplim,base_dcopf_dro:gen_dnlim} enforce the technical production limits of each controllable generator. 
Constraints \cref{base_dcopf_dro:lin_uplim,base_dcopf_dro:lin_dnlim} map the power injections and withdrawals of each resource and load to a resulting power flow using suitable linear maps \cite{dall2017chance}, e.g., obtained from the DC power flow approximation. 
Vectors $\bm{f}^{\rm RAM+}$ and $\bm{f}^{\rm RAM-}$ are the remaining available margins for each power transmission line, i.e., the distance of the power flow caused by the forecast injections and the upper and lower power flow limits.

Constraint \cref{base_dcopf_dro:joint_cc} enforces a joint chance constraint on the system response to uncertain injections $\bm{\xi}$. It ensures that the smallest probability, as measured by the distributions in ambiguity set $\set{A}$, for all constraints to hold is at least $1-\gamma$, where $\gamma$ is a risk parameter that is chosen by the system operator to be small. 
The first two constraints in \cref{base_dcopf_dro:joint_cc} ensure that control response does not exceed the available reserve and the last two constraints in \cref{base_dcopf_dro:joint_cc} ensure that resulting power flow changes do not exceed the remaining available margins on each power line. 
Finally, \cref{base_dcopf_dro:nonneg} lists all non-negative variables. 

For leaner notation in the following derivations, we collect all decision variables in a vector $\bm{x}$, define the space defined by constraints \cref{base_dcopf_dro:enerbal,base_dcopf_dro:resbal,base_dcopf_dro:gen_uplim,base_dcopf_dro:gen_dnlim,base_dcopf_dro:lin_uplim,base_dcopf_dro:lin_dnlim,base_dcopf_dro:nonneg} as $\set{F}$, and re-write \mbox{\cref{eq:base_dcopf_dro} as}
\begin{subequations}
\begin{align}
\min_{\bm{x}\in\set{F}}\ 
    &  \dprod{\bm{c}'}{\bm{x}} + \max_{\mathcal{Q}\in\set{A}} \mathbb{E}_{\mathcal{Q}}[\dprod{\bm{c}^{\rm A}}{-\bm{A}\bm{\xi}}] \label{compact_dcopf_dro:objective} \\
\text{s.t.}\quad
   & \inf_{\mathcal{Q}\in\set{A}} \mathcal{Q} \left\{\max_{k=1,...,K}[ \dprod{\bm{a}_{k}}{\bm{\xi}} + b_k] \le 0 \right\} \ge 1-\gamma \label{compact_dcopf_dro:joint_cc}
\end{align}%
\label{eq:compact_dcopf_dro}%
\end{subequations}%
where $\bm{a}_k$ and $b_k$ are, respectively, the $k$-th row and $k$-th entry of the 
$K\times D$ matrix and $K\times 1$ vector
\begin{equation*}
    \begin{bmatrix*}[c]
    -\bm{A} \\ \bm{A} \\ (\bm{B}^{\rm\bm{W}} - \bm{B}^{\rm G}\bm{A}) \\ -(\bm{B}^{\rm W} - \bm{B}^{\rm G}\bm{A})
    \end{bmatrix*}  \text{   and   } 
    \begin{bmatrix*}[c]
    -\bm{r}^+ \\ -\bm{r}^- \\ - \bm{f}^{\rm RAM+} \\ -\bm{f}^{\rm RAM-}
    \end{bmatrix*}.
\end{equation*}
and cost vector $\bm{c}$ combines $\bm{c}^{\rm E}$ and $\bm{c}^{\rm R}$. The formulations \cref{base_dcopf_dro:joint_cc} and \cref{compact_dcopf_dro:joint_cc} of the joint chance constraint are equivalent.

\subsection{Data and support}
\label{ssec:data_and_support}

To make inferences on the potential distributions of $\bm{\xi}$ as defined by ambiguity set $\set{A}$, the decision maker obtains datasets $\{\widehat{\set{X}}_j\}_{j=1}^D$.
These datasets contain information on each entry $\xi_j$ of $\bm{\xi}$ in the form of standardized samples $\{\hat{\xi}_{j1},...,\hat{\xi}_{jN'}\},\ j=1,...,D$ alongside their data-quality information $\{\epsilon_j\}_{j=1}^D$. 
{\color{black} For the remainder of the paper we assume that the $\epsilon_j$ are a confidence bound on the 1-Wasserstein distance ($p=1)$ and that the norm in \cref{eq:wasserstein_definition} is the 1-norm.
These assumptions will allow the formulation of a tractable linear program but are not restrictive in practice, because the data user can prescribe how the data owner should compute data quality.
}

We define support $\Xi$ as $\Xi_j = [\underline{\xi}_j, \overline{\xi}_j],\ j=1,...,D$ and assume $\underline{\xi}_j \le 0$ and $\overline{\xi}_j\ge 0$.
The support interval $[\underline{\xi}_j, \overline{\xi}_j]$ can be defined by the decision maker, e.g., by inferring it from the available data such that $\underline{\xi}_j = \min_i\{\widehat{\xi}_{ji}\}$, $\overline{\xi}_j = \max_i\{\widehat{\xi}_{ji}\}$.
For models of physical resources, a suitable definition of the support interval can be defined in terms of the physical limits $u_j^{\rm min}$, $u_j^{\rm max}$ of resource $j$. For a given forecast $u_j$ the technically feasible range of forecast errors is $[\kappa_j(u_j^{\rm min} - u_j), \kappa_j(u_j^{\rm max} - u_j)]$ \cite{esteban2021distributionally}, where $\kappa_j = [0,1]$ can be used as an additional scaling parameter. If $\kappa_j = 1$, the decision maker considers the entire technical range of the resource. If $\kappa_j = 0$ the decision maker effectively assumes that forecast $u_j$ will certainly be correct.
We follow \cite{esteban2021distributionally} and use this support definition for the remainder of this paper. 

\subsection{Objective reformulation}
\label{ssec:objective_reformulation}

We first reformulate the worst-case expectation in \cref{compact_dcopf_dro:objective} such that $\set{A}=\set{A}^{\rm MSW}$ using the results of Proposition~\ref{prop:mwsdro_sep_cost}, since the linear cost function $\dprod{\bm{c}^{\rm A}}{-\bm{A}\bm{\xi}}$ and support $\Xi$ are separable in $j$, which meets Condition~\ref{cond:separable_cost_and_support}
Applying Proposition~\ref{prop:mwsdro_sep_cost} to 
\cref{compact_dcopf_dro:objective} yields:
\allowdisplaybreaks
\begin{subequations}
\begin{align}
\min_{\bm{x}\in\set{F},  \lambda_j^{\rm co}, s_{ij}^{\rm co}}\  
    & \dprod{\bm{c}'}{\bm{x}} + 
    \sum_{j=1}^D \big(\lambda_j^{\rm co} \epsilon_j + \frac{1}{N_j}\sum_{i=1}^{N_j} s_{ji}^{\rm co}\big) \hspace{-3cm}&& \label{opf_wcexp:objective}\\
\text{s.t.}\
    & \forall j=1,...,D,\ i=1,...,N_j: \nonumber \\
 (\mu_{ji}^{\rm up}):\quad    
    & s_{ji}^{\rm co} \ge 
    -\sum_{g=1}^G c_g^{\rm A} \alpha_{gj} \overline{\xi}_j - \lambda_j^{\rm co} (\overline{\xi}_j - \widehat{\xi}_{ji}) \label{opf_wcexp:epi_u}\\  
(\mu_{ji}^{\rm lo}):\quad    
    & s_{ji}^{\rm co} \ge 
    -\sum_{g=1}^G c_g^{\rm A} \alpha_{gj} \underline{\xi}_j + \lambda_j^{\rm co} (\underline{\xi}_j - \widehat{\xi}_{ji}) \label{opf_wcexp:epi_l}\\
    & s_{ji}^{\rm co} \ge 
    -\sum_{g=1}^G c_g^{\rm A} \alpha_{gj}\widehat{\xi}_{ji}, \label{opf_wcexp:epi_av} \\
    & \lambda_j^{\rm co} \ge 0 \label{opf_wcexp:nonneg_lambda}
\end{align}%
\label{eq:opf_wcexp}%
\end{subequations}%
\allowdisplaybreaks[0]%
where $\alpha_{gj}$ is the entry in the $g$-th row and $j$-th column of matrix $\bm{A}$.
See Appendix~\ref{ax:wc_obj_derivation} for further derivation details. We introduce the superscript on $\lambda_j^{\rm co}$ to highlight its connection to the cost function reformulation.
Further, we note that the reformulation \cref{eq:opf_wcexp} of \cref{compact_dcopf_dro:objective} is \textit{exact}. 

\subsection{Chance constraint reformulation}
\label{ssec:chance_constraint_reformulation}

Without any assumptions on the distribution of $\bm{\xi}$, the joint chance constraint in  \cref{compact_dcopf_dro:joint_cc} cannot be solved directly.
Therefore, we use a common and convex inner approximation based on the conditional value-at-risk (CVaR) \cite{esteban2021distributionally,roveto2020co,dall2017chance,rockafellar2000optimization}. For \cref{compact_dcopf_dro:joint_cc} this approximation is given by:
\begin{equation}
    \max_{\mathcal{Q}\in\set{A}}\  \mathcal{Q}\cvar_{\gamma}\Big(\max_{k=1,...,K}[ \dprod{\bm{a}_{k}}{\xi} + b_k] \Big) \le 0,
\label{eq:cvar_base_constraint}
\end{equation}
where $\mathcal{Q}\cvar_{\gamma}$ is the CVaR at risk level $\gamma$ under distribution $\mathcal{Q}$. 
{\color{black} 
We note that the CVaR is an inner approximation of the VaR. Therefore, enforcing \cref{eq:cvar_base_constraint} implies that \cref{compact_dcopf_dro:joint_cc} holds \cite{rockafellar2000optimization}, thus ensuring the desired level of security $\gamma$.}
Enforcing \cref{eq:cvar_base_constraint} implies that \cref{compact_dcopf_dro:joint_cc} holds \cite{rockafellar2000optimization}, which ensures that the desired level of security $\gamma$ is met.
The CVaR formulation in \cref{eq:cvar_base_constraint} now allows an equivalent reformulation as
\begin{equation}
\begin{aligned}     
    & \inf_{\tau} \left\{ \tau + \frac{1}{\gamma}\max_{\mathcal{Q}\in\set{A}} \mathbb{E}_{\mathcal{Q}}\big[ \max_{k=1,...,K}\dprod{\bm{a}_{k}}{\bm{\xi}} + b_k - \tau \big]^+ \right\} \le 0\\
\Leftrightarrow\ 
    & \inf_{\tau} \left\{ \tau + \frac{1}{\gamma}\max_{\mathcal{Q}\in\set{A}} \mathbb{E}_{\mathcal{Q}}\big[ \max_{k=1,...,K+1}\dprod{\bm{a}_{k}'}{\bm{\xi}} + b_k'\big] \right\} \le 0 \\
\Leftrightarrow\ 
    &\left\{\begin{array}{l}
         \tau \le 0 \\
         \tau + \nu \le 0 \\
        \gamma \nu \ge \max_{\mathcal{Q}\in\set{A}} \mathbb{E}_{\mathcal{Q}}\big[ \max_{k=1,...,K+1}\dprod{\bm{a}_{k}'}{\bm{\xi}} + b_k'\big].
    \end{array}\right.
\end{aligned}%
\label{eq:joint_cc_cvar_approx}%
\end{equation}%
The initial formulation in the first line of \cref{eq:joint_cc_cvar_approx} follows from the dual form of the CVaR \cite{rockafellar2000optimization}.
The second line of \cref{eq:joint_cc_cvar_approx} introduces $\bm{a}'_k = \bm{a}'_k(\bm{x})$ and $b'_k = b'_k(\bm{x},\tau)$ such that $\bm{a}'_k = \bm{a}_k,\ k=1,...,K$, $\bm{a}'_{K+1} = \bm{0}_D$, $b'_k = b_k-\tau,\ k=1,...,K$, $b'_{K+1} = 0$, where $\bm{0}_D$ denotes a vector of zeros of length $D$.  
The third line of \cref{eq:joint_cc_cvar_approx} allows for the limiting case $\gamma=0$.

We can now reformulate the inner worst-case expectation in \cref{eq:joint_cc_cvar_approx} such that $\set{A}=\set{A}^{\rm MSW}$ using the results of Proposition~\ref{prop:mwsdro_stand_data}. We apply this proposition here because the function inside the expectation operator of \cref{eq:joint_cc_cvar_approx} cannot be separated into its $\xi_j$ components due to the inner maximization operator.
Using Proposition~\ref{prop:mwsdro_stand_data}, the worst-case expectation in the last line of \cref{eq:joint_cc_cvar_approx} becomes:
\allowdisplaybreaks
\begin{subequations}
\begin{align}
\inf_{\lambda_j^{\rm cc}\ge0}\ 
    & \sum_{j=1}^D \lambda_j^{\rm cc} \epsilon_j + \frac{1}{N'}\sum_{i=1}^{N'} s_{i}^{\rm cc} \\
\text{s.t.} \quad 
    & s_{i}^{\rm cc}\! \ge \max_{\bm{\xi}\in\Xi} \Big( \!\max_{k=1,...,K+1}\!(\dprod{\bm{a}_k'}{\bm{\xi}} \!+\! b_k')\!-\!\! \sum_{j=1}^D\lambda_j^{\rm cc} |\xi_j\! - \widehat{\xi}_{ji'}|\!\Big), \nonumber \\
    & \hspace{4.5cm}    i=1,...,N',
    \label{seq2:cvar_reform1}
\end{align}%
\label{eq:joint_cc_cvar_reform1}%
\end{subequations}%
\allowdisplaybreaks[0]%
where we introduce the superscript on $\lambda_j^{\rm cc}$ to highlight its connection to the chance-constraint reformulation. 
We now find a tractable reformulation of \cref{seq2:cvar_reform1} as follows. 
First, we move the inner maximization outwards, enforce it through an epigraph formulation, and rewrite the inner dot product:
\begin{align}
    & s_{i}^{\rm cc}\! \ge b_k' + \max_{\bm{\xi}\in\Xi} \sum_{j=1}^D \big(a_{kj}'\xi_j - \lambda_j^{cc}|\xi_j\! - \widehat{\xi}_{ji'}|)\big), \nonumber \\
    & \hspace{2cm}  i\!=\!1,...,N',\ k\!=\!1,...,K\!\!+\!1. \label{eq:joint_cc_cvar_reform2}%
\end{align}%
The remaining inner supremum in \cref{eq:joint_cc_cvar_reform2} can be resolved by realizing that the supremum is obtained when each term in the sum is maximized and that either $\overline{\xi}_j$, $\underline{\xi}_j$, or $\widehat{\xi}_{ij}$ will be the maximizer.
This allows the equivalent reformulation of \cref{eq:joint_cc_cvar_reform2}:
\begin{subequations}
\begin{align} 
(\eta_{ik}):\quad& s_i^{\rm cc} \ge b_k' + \sum_{j=1}^D s_{jik},\ i=1,...,N',\ k\!=\!1,...,K\!\!+\!1\label{cc_reform:main}\\
&\hspace{-1.5cm}j=1,...,D,\ i=1,...,N',\ k=1,...,K+1:  \nonumber \\
(\rho_{jik}^{\rm up}):\quad &s_{jik} \ge a_{kj}'\overline{\xi}_j\!\! -\!\! \lambda_{j}^{\rm cc} (\overline{\xi}_j - \widehat{\xi}_{ji}), \label{cc_reform:aux_up} \\ 
(\rho_{jik}^{\rm lo}):\quad& s_{jik} \ge a_{kj}'\underline{\xi}_j\! \!+\! \!\lambda_{j}^{\rm cc} (\underline{\xi}_j - \widehat{\xi}_{ji}),  \label{cc_reform:aux_lo}\\
(\rho_{jik}^{\rm av}):\quad& s_{jik} \ge a_{kj}'\widehat{\xi}_{ji}, \label{cc_reform:aux_av}
\end{align}%
\label{eq:joint_cc_cvar_final}%
\end{subequations}%
where we introduce $s_{jik}$ as auxiliary variables.

\subsection{Complete formulation}
The final model formulation is given by:
\allowdisplaybreaks
\begin{subequations}
\begin{align}
&\min_{\substack{\bm{x}, \lambda_j^{\rm co}, s_{ij}^{\rm co}\\ \tau\le0, \nu, \lambda_{j}^{\rm cc}, s_i^{\rm cc} }
}\  
     \dprod{\bm{c}'}{\bm{x}} + 
    \sum_{j=1}^D \big(\lambda_j^{\rm co} \epsilon_j + \frac{1}{N'}\sum_{i=1}^{N'} s_{ji}^{\rm co}\big) \label{complete_msdro_opf:objective}\\
\text{s.t.}\ &\text{\cref{base_dcopf_dro:enerbal,base_dcopf_dro:resbal,base_dcopf_dro:gen_uplim,base_dcopf_dro:gen_dnlim,base_dcopf_dro:lin_uplim,base_dcopf_dro:lin_dnlim,base_dcopf_dro:nonneg} [Deterministic constraints]} \nonumber \\
&\text{\cref{opf_wcexp:epi_u,opf_wcexp:epi_l,opf_wcexp:epi_av,opf_wcexp:nonneg_lambda} [Worst-case exp. cost aux. constraints]} \nonumber \\
 &\quad \tau + \nu \le 0 \label{complete_msdro_opf:cvar_const1} \\ 
(\phi): &\quad \sum_{j=1}^{D} \epsilon_j \lambda_j^{\rm cc} 
 + \frac{1}{N'} \sum_{i=1}^{N'} s_i^{\rm cc} \le \gamma \nu \label{complete_msdro_opf:cvar_const2}\\
 & \text{\cref{cc_reform:main,cc_reform:aux_up,cc_reform:aux_lo,cc_reform:aux_av} [CVaR aux. constraints]}. \nonumber
\end{align}%
\label{eq:complete_msdro_opf}%
\end{subequations}%
\allowdisplaybreaks[0]%

{\color{black}
\begin{remark}
The formulations in this section assume that the $\epsilon_j$ are reported accurately, i.e., as discussed in Section~\ref{sec:privacy_and_dist_ambiguity}.
In practice, it is possible that the data owner reports data quality incorrectly. 
In that case, additional techniques are required to detect such inaccuracies, to correct them, and to provide mechanisms that incentivize data owners to report data quality to the best of their ability.
However, the design of such techniques is beyond the scope of this paper and subject of further research.
\end{remark}
}

\section{Data valuation}
\label{sec:analyis}

{\color{black}%
We now focus on the relationship between the data quality information encoded in $\epsilon_j$ and the \textit{value} of a corresponding dataset $\widehat{\set{X}}_j$ for a decision-making problem. 
Because our formulations derived in Theorem~\ref{th:multi-source-wasserstein-dro} and Propositions~\ref{prop:mwsdro_sep_cost} and \ref{prop:mwsdro_stand_data} allow to internalize the individual data quality of each dataset, we can now explicitly formulate the \textit{marginal value of data quality}, i.e., the marginal improvement of the objective as a function of $\epsilon_j$.
Below we provide a detailed analysis of this concept for the OPF problem derived in Section~\ref{sec:ms_opf}.
}

\subsection{Marginal value of data quality}
\label{ssec:marginal_data_value_analysis}

Variables $\lambda_j^{\rm co}$ and $\lambda_j^{\rm cc}$ allow an interpretation as the \textit{marginal value of data quality} in the optimal solution of \cref{eq:complete_msdro_opf}. 
Consider $\mathcal{L}$ as the Lagrangian of \cref{eq:complete_msdro_opf} and denote $\phi$ as the dual multiplier of \cref{complete_msdro_opf:cvar_const2}.
For a given (primal and dual) optimal solution of \cref{complete_msdro_opf:cvar_const2}, we can employ the Envelope Theorem to determine the sensitivity of the optimal objective value to the data quality as 
\begin{equation}
    \frac{\partial \mathcal{L}}{\partial \epsilon_j} = \lambda_j^{\rm co} + \phi \lambda_j^{\rm cc}.
    \label{eq:L_over_epsilon}
\end{equation}
This term is composed of the direct impact of data quality on the uncertain part of the objective, captured by $\lambda_j^{\rm co}$, and the indirect impact of data quality on the objective through the chance constraint, captured by $\phi \lambda_j^{\rm cc}$.
Here, $\phi$ allows an interpretation as the marginal cost of enforcing all constraints within the joint chance constraint as per: 
\begin{equation}
\begin{aligned}
&\frac{\partial \mathcal{L}}{\partial s_i^{\rm cc}} = \frac{\phi}{N'} - \sum_{k=1}^{K+1} \eta_{ik} \stackrel{\text{(A)}}{=} 0\ 
\Rightarrow \  \phi =  N'\sum_{k=1}^{K+1} \eta_{ik}, \quad \forall i,
\end{aligned}
\end{equation}
where $\eta_{ik}$ is the dual multiplier of \cref{cc_reform:main} and relation (A) follows from the optimality conditions of \cref{eq:complete_msdro_opf}. In other words, $\phi$ translates the per-feature sensitivity $\lambda_j^{\rm cc}$ of all constraints to data quality $\epsilon_j$ into a cost value.

Variables $\lambda_j^{\rm co}$ and $\lambda_j^{\rm cc}$ also allow a physical interpretation. Consider $\lambda_j^{\rm co}$ and the reformulation of the objective function in \cref{eq:opf_wcexp}.
In its second term, objective \cref{opf_wcexp:objective} needs to balance smaller values of $\lambda_j^{\rm co}$ with smaller values of $s^{\rm co}_{ji}$. 
The lowest possible value of $\lambda_j^{\rm co}$ is $\lambda_j^{\rm co}=0$. 
In this case, the second terms on the right-hand side of \cref{opf_wcexp:epi_l,opf_wcexp:epi_u} disappear and the epigraph formulation in  \cref{opf_wcexp:epi_l,opf_wcexp:epi_u,opf_wcexp:epi_av} sets $s^{\rm co}_{ji} = -\sum_{g=1}^G c_g^{\rm A}\alpha_{gj}\underline{\xi}_j,\ \forall i$, i.e., the \textit{robust solution} given support $\Xi_j$. 
Defaulting to the robust solution for feature $j$ means that available data points $\{\widehat{\xi}_{ji}\}_{i=1}^{N}$ are ignored due to insufficient data quality, i.e., $\epsilon_j$ was too large to offer the decision-maker any information beyond the known support. 
Hence, the value of the data is zero. 
On the other hand, the largest possible value of $\lambda_j^{\rm co}$ in \cref{eq:opf_wcexp} is $\lambda_j^{\rm co} = \sum_{g=1}^G c_g^{\rm A} \alpha_{gj}$.
In this case constraints \cref{opf_wcexp:epi_l,opf_wcexp:epi_u} are equal to \cref{opf_wcexp:epi_av} and the epigraph formulation in  \cref{opf_wcexp:epi_l,opf_wcexp:epi_u,opf_wcexp:epi_av} sets $s^{\rm co}_{ji} = -\sum_{g=1}^G c_g^{\rm A}\alpha_{gj}\widehat{\xi}_j,\ \forall i$, which indicates that the marginal value of data quality equals the marginal cost of reserve activation.
This result offers a complementary physical interpretation of $\epsilon_j$, which we discuss in the following Section~\ref{ssec:data_quality_interpretaion}. 

\begin{remark}
An analogous reasoning is possible for $\lambda_j^{\rm cc}$ showing that  $\lambda_j^{\rm cc}$ is either zero, indicating that the joint chance constraint defaults to the robust solution for feature $j$, or $\exists k$ such that $\lambda_j^{\rm cc} = |a'_{kj}|$. 
We omit the formal derivation of this result here for brevity of this article.
\end{remark}

\subsection{Analysis of data quality}
\label{ssec:data_quality_interpretaion}
We again focus on the analysis of $\lambda_j^{\rm co}$ due to its more accessible interpretation.
Consider $\underline{\xi}_j$, which, as discussed above, is the worst-case (robust) corner of support $\Xi$ in $j$.
\begin{proposition}
\label{prop:epsilon_analysis}
Let $\lambda_j^{\rm co,*}$ be the value of $\lambda_j^{\rm co}$ in the optimal solution of \cref{eq:complete_msdro_opf}. It holds that:
\begin{equation*}
    \lambda_j^{\rm co, *} = \begin{cases}
        0, &\text{if }  \epsilon_j \ge \frac{1}{N'} \sum_{i=1}^{N'}(\widehat{\xi}_{ji} - \underline{\xi}_j) \\
        \sum_{g=1}^G c_g^{\rm A} \alpha_{gj}, &\text{else.}
    \end{cases}
\end{equation*}
\end{proposition}
\myproofstart
    See Appendix~\ref{ax:proof_epsilon_analysis}.
\myproofend

Proposition~\ref{prop:epsilon_analysis} offers two valuable insights.
First, it states that $\lambda_j^{\rm co, *}$ is either zero or $\sum_{g=1}^G c_g^{\rm A} \alpha_{gj}$.
Hence, if data-quality is sufficient to allow a non-zero $\lambda_j^{\rm co}$, then the worst-case expected cost caused by feature $j$ is always equal to $\sum_{g=1}^G c_g^{\rm A}\alpha_{gj}(\epsilon_j - \frac{1}{N'}\sum_{i=1}^N\widehat{\xi}_{ji})$.
This implies a physical interpretation of $\epsilon_j$ as the \textit{additional expected amount of reserve activation beyond the sample average}.
Second, Proposition~\ref{prop:epsilon_analysis} states that $\lambda_j^{\rm co, *}$ is zero whenever $\epsilon_j$ is larger than the average distance of the data from the worst-case boundary of the support. 
This allows the decision-maker to assess data usefulness \textit{offline}, i.e., before solving the actual decision-making problem.

\subsection{Towards data remuneration}
\label{ssec:data_remuneration}

The interpretation of $\lambda_j^{\rm co}$ and $\lambda_j^{\rm cc}$ as a marginal cost aligns with the established theory and many proposed approaches for pricing and remuneration schemes in power systems \cite{kazempour2018stochastic,fang2019introducing,mieth2019distribution,kuang2018pricing}. 
Our proposed formulation allows for similar formulations that can capture the value of an uncertain resource given the data used to inform its forecast and uncertainty. 
For example, in \cite{fang2019introducing} an uncertain resource is remunerated based on the marginal value of its forecast and a direct relationship between the forecast and the resulting uncertainty (i.e., variance of the forecast) is implied. In this way, the value of the uncertain resource can be decomposed into a ``firm'' value contribution related to the forecast and a value ``deduction'' related to the uncertainty cost the resource introduces into the system.
Using the results from Sections~\ref{ssec:marginal_data_value_analysis} and \ref{ssec:data_quality_interpretaion}, our approach achieves a parallel result without imposing a functional relationship between the forecast and its uncertainty. 

Consider forecast value $u_j$ and its marginal value given by
\allowdisplaybreaks
\begin{align}
\pi_j^{\rm F} \coloneqq  \frac{\partial \mathcal{L}}{\partial u_j} \!=\!
& \underbrace{\pi + \sum_{l=1}^L B_{lj}^{\rm W}(\beta_{l}^{up} - \beta_{l}^{lo})}_{\text{(\ref{eq:L_over_u}a): locational marginal price}} \nonumber \\
&\hspace{-2.2cm}- \!\underbrace{\kappa_j\! \sum_{i=1}^{N'}\!\Big(\mu_{ji}^{\rm up}\!\big(\!\sum_{g=1}^G\!c_{g}^A\alpha_{gj}\! +\! \lambda_j^{\rm co}\big)\! +\! \mu_{ji}^{\rm lo}\!\big(\!\sum_{g=1}^G\!c_{g}^A\alpha_{gj}\! -\! \lambda_j^{\rm co}\big)\!\Big)}_{\text{(\ref{eq:L_over_u}b): contribution to expected balancing cost}} \nonumber \\
&\hspace{-2.2cm} -\!\underbrace{\kappa_j \!\sum_{i=1}^{N'} \sum_{k=1}^K \Big(\rho_{jik}^{\rm up} (\lambda_j^{\rm cc} - a_{kj}') - \rho_{jik}^{\rm lo} (\lambda_j^{\rm cc} + a_{kj}' )\!\Big).}_{\text{(\ref{eq:L_over_u}c): contribution to system cost from reserves}}\label{eq:L_over_u}
\end{align}%
\allowdisplaybreaks[0]%
Term (\ref{eq:L_over_u}a) is the locational marginal price (LMP) of the bus to which resource $j$ is connected and reflects the value of the forecast power injection $u_j$ without considering its uncertainty. 
Similar to the result in \cite{fang2019introducing}, terms  (\ref{eq:L_over_u}b) and (\ref{eq:L_over_u}c) reduce this value based on how the resource contributes to the cost of providing reserve for its uncertainty. If $\kappa_j=0$, i.e., the decision maker assumes that the forecast is certainly true (see Section~\ref{ssec:data_and_support}), the value of the resources collapses to its LMP.
However, following results from Sections~\ref{ssec:marginal_data_value_analysis} and \ref{ssec:data_quality_interpretaion}, terms (\ref{eq:L_over_u}b) and (\ref{eq:L_over_u}c) are mainly determined by the support definition. They approach zero when $\lambda_j^{\rm co} = \sum_{g=1}^Gc_{g}^A\alpha_{gj}$, $\lambda_j^{\rm cc} = |a_{kj}'|$ and are equal to the robust solution if $\lambda_j^{\rm co} = 0$, $\lambda_j^{\rm cc} = 0$.
We can correct the resulting missing value reduction by adding $\epsilon_j \frac{\partial \mathcal{L}}{\partial \epsilon_j}$ (see \cref{eq:L_over_epsilon}) to the resource valuation, due to the interpretation of $\epsilon_j$ as the critical additional amount of uncertainty (see Section~\ref{ssec:data_quality_interpretaion}) and the interpretation of $\lambda_j^{\rm co}$, $\lambda_j^{\rm cc}$ capturing the marginal cost of uncertainty  (Section~\ref{ssec:marginal_data_value_analysis}).
This offers a pathway to remunerating data-informed \mbox{resources as} 
\begin{equation}
    u_j \pi_j^{\rm F}(u_j) - \epsilon_j(\lambda_j^{\rm co} + \phi\lambda_j^{\rm cc}) \eqqcolon  u_j\pi_j^{\rm F} - \epsilon_j\pi_j^{\rm D}.
\label{eq:complete_remuneration}
\end{equation}
An immediate result of \cref{eq:complete_remuneration} is that it encourages a low $\epsilon_j$, i.e., high data quality.
We discuss \cref{eq:complete_remuneration} in the case study below, but note that further analysis of properties and potential market designs are out of scope here and subject to future work.

\section{Case study}

We conduct illustrative numerical experiments for the multi-source data-driven DRO OPF \cref{eq:complete_msdro_opf} using a modified version of the IEEE 5-bus test system. 
We choose this small-scale example to accommodate a more accessible analysis of the proposed formulation. We discuss scalability and larger-scale applications in Sections~\ref{ssec:scalability} and \ref{sec:conlcusion}. 

\subsection{Data and implementation}

\begin{figure}
    \centering
    \includegraphics[width=0.8\linewidth]{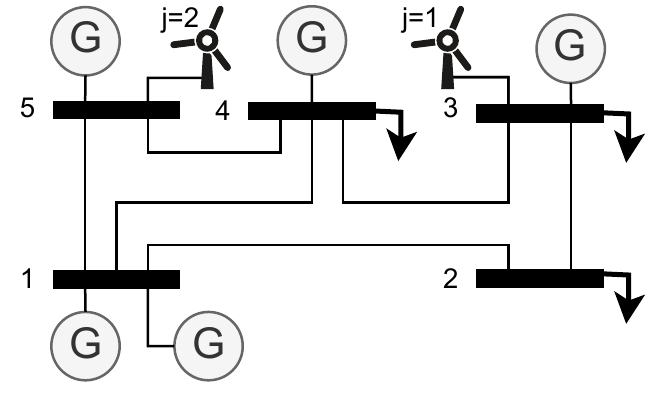}
    \caption{Schematic of the 5 bus test system.}
    \label{fig:5bus_system}
\end{figure}

We use the ``case5'' dataset from MATPOWER \cite{matpowercase5} with the topology shown in Fig.~\ref{fig:5bus_system}. We enforce thermal limits on all lines with $\bm{f}^{\rm max}=[3.2, 1.52, 1.76, 0.8, 0.8, 1.92]\tran{\rm p.u.}$ and add two uncertain resources with forecasts $\bm{u}=[1,1.5]{\rm p.u.}$ at buses 3 and 4 as shown in Fig.~\ref{fig:5bus_system} as wind turbines.
We set the maximum capacity $u_j^{\rm max}$ of both uncertain resources to \unit[2]{p.u.} and $\kappa_j=\kappa=0.6,\ j=1,2$ such that the resulting maximum forecast error support is $\underline{\bm{\xi}} = -\kappa\bm{u} = [-0.6,-0.9]\tran{\rm p.u.}$ and $\overline{\bm{\xi}} = \kappa(\bm{u}^{\rm max}-\bm{u})=[0.6,0.3]\tran{\rm p.u.}$
The generator cost vectors are set to $\bm{c}^{\rm E}=[14,15,30,40,10]\tran$, $\bm{c}^{\rm A}=10\bm{c}^{\rm R}=[80,80,15,30,80]\tran$, where $\bm{c}^{\rm E}$ is part of the original dataset and $\bm{c}^{\rm R}$, $\bm{c}^{\rm A}$ have been chosen to prioritize generators with higher energy cost (e.g., gas, oil) for balancing reserves. 

For each wind farm $j=\{1,2\}$ we generate a set $\widehat{\set{X}}_j$ of $N_j = N' = 20,\ \forall j$ data samples to inform forecast errors by sampling data from a normal distribution with mean $u_j$ (i.e., the forecast) and standard deviation $S_j = 0.15 u_j$ \cite{dvorkin2015uncertainty}. 
Further, the distribution is truncated to reflect the support interval $[\underline{\xi}_j, \overline{\xi}_j]$ for each resource $j$. We note that for the purpose of illustrating the various components of the proposed model and their relationship, we do not explicitly model raw or private datasets $\set{Y}_j$ or $\set{X}_j$ (see Section~\ref{sec:privacy_and_dist_ambiguity}) here, but assign data quality assumptions $\epsilon_j$ to all $\widehat{\set{X}}_j$ \textit{ex-post}.
We set $\gamma = 5\%$.

Our implementation in the Julia language using JuMP \cite{lubin2023jump} with the Gurobi solver is available open source \cite{git_msdro_opf}.
All computations have been performed on a standard PC workstation with \unit[16]{GB} memory and an Intel i5 processor. 

\begin{remark}
While the formulation in \cref{eq:complete_msdro_opf} can be solved directly with off-the-shelf solvers, our implementation uses a simple extra step to improve the tightness of the CVaR approximation of the joint chance-constraint (see Section~\ref{ssec:chance_constraint_reformulation}).
If the model decides that a set of generators $\set{G}'$ should not participate in balancing, i.e., $\alpha_{gj}=0,\ g\in\set{G}',\ \forall j$, then their balancing contribution $\dprod{\bm{\alpha}_g}{\bm{\xi}} = 0,\ \forall \bm{\xi}$, implying $r^+ = r^- =0$.
We can thus improve the CVaR approximation by re-running the model, setting $r_g^+ = r_g^- = 0,\ g\in\set{G}'$, and moving the respective reserve constraints outside of the CVaR constraint. 
This correction improves the accuracy of the CVaR approximation similar to the constraint-scaling methods discussed in \cite{ordoudis2021energy}.

\end{remark}

\subsection{Results and analysis}
\label{ssec:results_and_analysis}

We run our numerical experiment 16 times for all combinations of data quality values $\epsilon_j\in\{1.0,0.1,0.005,0.001\}$.
Tables~\ref{tab:objective_values} and \ref{tab:lambda_values} show the resulting values for the expected system cost and $\lambda_j^{\rm co}$, $\lambda_j^{\rm cc}$, respectively.
For $(\epsilon_1, \epsilon_2)=(1.0,1.0)$ we observe that $\lambda_j^{\rm co}=\lambda_j^{\rm cc}=0$, i.e., data quality is insufficient and generation and transmission reserves $r^+$, $r^-$, $f^{\rm RAM+}$, $f^{\rm RAM-}$ are sized to provide a robust solution within support $\Xi = [\underline{\xi}_1, \overline{\xi}_1] \times [\underline{\xi}_2, \overline{\xi}_2]$ leading to the highest expected system cost (see Table~\ref{tab:objective_values}). With smaller $\epsilon_j$, indicating smaller required ambiguity sets around the empirical data distributions, the operator can put higher trust in the data and reduce the uncertainty space around the forecast values. This leads to reduced reserves, lower system cost, and an increase of  \mbox{$\lambda_j^{\rm co}$, $\lambda_j^{\rm cc}$}. 

\begin{table}[]
\caption{Objective values [\$] for various $\epsilon_j$}
\label{tab:objective_values}
\renewcommand{\arraystretch}{1.1}
    \centering
\begin{tabular}{c|cccc}
 \diagbox{$\epsilon_1$}{$\epsilon_2$}   & \textbf{1.0} & \textbf{0.1} & \textbf{0.005} & \textbf{0.001}\\
    \hline
    \textbf{1.0} & 24241.6 & 19293.9 & 18339.2 & 17811.0\\
    \textbf{0.1} & 23491.6 & 18543.9 & 17589.2 & 17061.0\\
    \textbf{0.005} & 23349.1 & 18401.4 & 17446.7 & 16918.5\\
    \textbf{0.001} & 23343.1 & 18395.4 & 17451.2 & 16912.5\\
\hline
\end{tabular}
\end{table}

\begin{table}[]
\renewcommand{\arraystretch}{1.1}
    \caption{Results for $\lambda_j^{\rm co}$ and $\lambda_j^{\rm cc}$ for various $\epsilon_j$}
    \label{tab:lambda_values}
    \centering
\begin{tabular}{cc|cc|cc}
   $\epsilon_1$ & $\epsilon_2$ & $\lambda_1^{\rm co}$ & $\lambda_2^{\rm co}$ & $\lambda_1^{\rm cc}$ & $\lambda_2^{\rm cc}$ \\
   \hline
   \textbf{1.0} & \textbf{1.0} &0 & 0 & 0.000 & 0.000 \\ 
   \textbf{1.0} & \textbf{0.1} &0 & 6901 & 0.000 & 0.000 \\ 
   \textbf{1.0} & \textbf{0.005} &0 & 5965 & 0.000 & 0.313 \\ 
   \textbf{1.0} & \textbf{0.001} &0 & 8000 & 0.000 & 1.000 \\ 
   \textbf{0.1} & \textbf{1.0} &1500 & 0 & 0.000 & 0.000 \\ 
   \textbf{0.1} & \textbf{0.1} &1500 & 6901 & 0.000 & 0.000 \\ 
   \textbf{0.1} & \textbf{0.005} &1500 & 5965 & 0.000 & 0.313 \\ 
   \textbf{0.1} & \textbf{0.001} &1500 & 8000 & 0.000 & 1.000 \\ 
   \textbf{0.005} & \textbf{1.0} &1500 & 0 & 0.000 & 0.000 \\ 
   \textbf{0.005} & \textbf{0.1} &1500 & 6901 & 0.000 & 0.000 \\ 
   \textbf{0.005} & \textbf{0.005} &1500 & 5965 & 0.000 & 0.313 \\ 
   \textbf{0.005} & \textbf{0.001} &1500 & 8000 & 0.000 & 1.000 \\ 
   \textbf{0.001} & \textbf{1.0} &1500 & 0 & 0.000 & 0.000 \\ 
   \textbf{0.001} & \textbf{0.1} &1500 & 6901 & 0.000 & 0.000 \\ 
   \textbf{0.001} & \textbf{0.005} &1665 & 5606 & 0.110 & 0.250 \\ 
   \textbf{0.001} & \textbf{0.001} &1500 & 8000 & 0.000 & 1.000 \\ 
\hline
\end{tabular}
\end{table}

We see in Table~\ref{tab:lambda_values} that the relationship between $\epsilon_j$ and $\lambda_j^{\rm co}$, $\lambda_j^{\rm cc}$ is not necessarily monotonic, e.g., in the transition from $(\epsilon_1, \epsilon_2)=(1.0,0.1)$ to $(\epsilon_1, \epsilon_2)=(1.0,0.005)$. 
This can be ascribed to the relation between data value to the physical decision-making problem, as discussed in Section~\ref{sec:analyis}.
We see in Table~\ref{tab:generator_production_and_reserve} that, with increasing data quality from resource $j=2$, the system is able to reduce overall reserves ($\sum_{g=1}^Gr_g = 1.5{\rm p.u.}$ for $(\epsilon_1, \epsilon_2)=(1.0,0.1)$ and $\sum_{g=1}^Gr_g = 1.22{\rm p.u.}$ for $(\epsilon_1, \epsilon_2)=(1.0,0.005)$) and shift balancing participation from the more expensive reserve provider $g=5$ to the cheaper $g=3$. As a result $\lambda_2^{\rm co}$ reduces, as per the interpretation provided in Section~\ref{ssec:marginal_data_value_analysis}.  
Table~\ref{tab:cost_components} supports this result by showing the equivalency of $\lambda_j^{\rm co}$ and $\sum_{g=1}^Gc_g^{\rm A}\alpha_{gj} = (\bm{c}^{\rm A})\tran \bm{A}_j$ whenever $\lambda_j^{\rm co}>0$. 
Simultaneously, $\lambda_2^{\rm cc}$ increases, showing that data quality is now sufficient to depart from the robust solution in the joint chance constraint, and highlighting the dependency between the two data value components $\lambda_j^{\rm co}$ and $\phi\lambda_j^{\rm cc}$.
For further reduction of $\epsilon_2$ to $\epsilon_2=0.001$, however, the system can further reduce reserve, but prefers to fully allocate balancing reserve for uncertain resource $j=2$ to the more expensive reserve provider $g=5$ (i.e., $\alpha_{52}=1.0$, not shown in Table~\ref{tab:generator_production_and_reserve}) leading to a lower system cost (see Table~\ref{tab:objective_values}) and a higher $\lambda_2^{\rm co} = \sum_{g=1}^Gc_g^{\rm A}\alpha_{gj}$ (see Tables~\ref{tab:lambda_values} and \ref{tab:cost_components}).
The results in Table~\ref{tab:lambda_values} further highlight that the marginal value of a dataset also depends on the quality of other datasets. We can observe this for $\epsilon_1=0.001$ at the bottom of Table~\ref{tab:lambda_values} where decreasing $\epsilon_2$ changes both $\lambda_1^{\rm co}$ and $\lambda_1^{\rm cc}$.

\begin{table}[]
\setlength{\tabcolsep}{4pt}
\renewcommand{\arraystretch}{1.15}
    \caption{Generator production and reserve}
    \label{tab:generator_production_and_reserve}
    \centering
    \begin{tabular}{c|cc|cccc|cccc}
         \multicolumn{3}{c|}{} & \multicolumn{4}{c|}{$(\epsilon_1, \epsilon_2) = (1.0, 0.1)$} & \multicolumn{4}{c}{$(\epsilon_1, \epsilon_2) = (1.0, 0.005)$} \\
        \hline
         & $c_g^{\rm E}$ & $c_g^{\rm A}$ & $p_g$ & $r_g$ & $\alpha_{g1}$ & $\alpha_{g2}$ & $p_g$ & $r_g$ & $\alpha_{g1}$ & $\alpha_{g2}$ \\
         g & \multicolumn{2}{c|}{[$\nicefrac{100\$}{\rm p.u.}$]} & [p.u.] & [p.u.] & [1]& [1]& [p.u.] & [p.u.] & [1]& [1] \\
        \hline 
\textbf{1} & 14 & 80 & 0.40\textsuperscript{*} & 0.00 & 0.00 & 0.00 & 0.40\textsuperscript{*} & 0.00 & 0.00 & 0.00\\
\textbf{2} & 15 & 80 & 1.70\textsuperscript{*} & 0.00 & 0.00 & 0.00 & 1.70\textsuperscript{*} & 0.00 & 0.00 & 0.00\\
\textbf{3} & 30 & 15 & 2.34 & 0.75 & 1.00 & 0.17 & 2.31 & 0.74 & 1.00 & 0.31\\
\textbf{4} & 40 & 30 & 1.03 & 0.00 & 0.00 & 0.00 & 1.10 & 0.00 & 0.00 & 0.00\\
\textbf{5} & 10 & 80 & 2.03 & 0.75 & 0.00 & 0.83 & 2.00 & 0.48 & 0.00 & 0.69\\
    \hline
    \end{tabular}\\
    \flushleft
   * Generator operates at $p_g^{\rm max}$
   ~~~~~~~~~~~~~~~~~~~~~~~~~~~~~~~
    $r_g = r_g^+ + r_g^-$
\end{table}

The last 8 columns in Table~\ref{tab:cost_components} itemize results related to the cost components studied in Eqs.~\cref{eq:L_over_u,eq:complete_remuneration}. When $\lambda_j^{\rm co}=\lambda_j^{\rm cc}=0,\ \forall j$, as in the case for $(\epsilon_1, \epsilon_2)=(1.0,1.0)$, the value of a data-informed uncertain resource is reduced only by the terms (\ref{eq:L_over_u}b), columns 5 and 6 in Table~\ref{tab:cost_components}, and (\ref{eq:L_over_u}c), columns 9 and 10 in Table~\ref{tab:cost_components}. 
As discussed in Section~\ref{ssec:data_remuneration}, for each $j$ (\ref{eq:L_over_u}b) is zero when $\lambda_j^{co}>0$ and the uncertainty cost (here for reserve activation) is captured by $\epsilon_j\lambda_j^{\rm co}$. We make the same observation for (\ref{eq:L_over_u}c) and $\epsilon_j\phi\lambda_j^{\rm cc}$. However, because the joint chance-constraint \cref{base_dcopf_dro:joint_cc} can not be separated in $j$, i.e., it depends on the combined impact of all uncertain resources, the robust solution of $j=2$ dominates the robust solution of $j=1$ leading its zero cost contribution from (\ref{eq:L_over_u}c).

\begin{table}[]
\setlength{\tabcolsep}{3pt}
\renewcommand{\arraystretch}{1.15}
\caption{Cost components for various $\epsilon_j$}
\label{tab:cost_components}
    \centering
    \begin{tabular}{cc|cc|cc|cc|cc|cc}
    &  & \multicolumn{2}{c|}{$(\bm{c}^{\rm A})\tran \bm{A}_j$} & \multicolumn{2}{c|}{$u_j$(\ref{eq:L_over_u}b)}  & \multicolumn{2}{c|}{$\epsilon_j\lambda_j^{\rm co}$} & \multicolumn{2}{c|}{$u_j$(\ref{eq:L_over_u}c)} & \multicolumn{2}{c}{$\epsilon_j\lambda_j^{\rm cc}\phi$} \\
    $\epsilon_1$ & $\epsilon_2$ & j=1 & j=2 & 1 & 2 & 1 & 2 & 1 & 2 & 1 & 2 \\
    \hline
\textbf{1.0} & \textbf{1.0} & 1500 & 5301 & 900 & 4771 & 0 & 0 & 0 & 1457 & 0 & 0\\
\textbf{1.0} & \textbf{0.1} & 1500 & 6901 & 900 & 0 & 0 & 690 & 0 & 304 & 0 & 0\\
\textbf{1.0} & \textbf{0.005} & 1500 & 5965 & 900 & 0 & 0 & 30 & 0 & 0 & 0 & 304\\
\textbf{1.0} & \textbf{0.001} & 1500 & 8000 & 900 & 0 & 0 & 8 & 0 & 0 & 0 & 163\\
\textbf{0.1} & \textbf{0.001} & 1500 & 8000 & 0 & 0 & 150 & 8 & 0 & 0 & 0 & 163\\
\textbf{0.001} & \textbf{0.005} & 1665 & 5606 & 0 & 0 & 2 & 28 & 0 & 0 & 22 & 256\\
\textbf{0.001} & \textbf{0.001} & 1500 & 8000 & 0 & 0 & 2 & 8 & 0 & 0 & 0 & 163\\
\hline
    \end{tabular}\\
\flushleft
 Cols. 3 and 4 in [\$/p.u.], rest in [\$].
\end{table}

{\color{black}
Table~\ref{tab:oos_joint_cc} shows an out-of-sample analysis of the empirical probability of one of the constraints inside the joint chance constraint \cref{base_dcopf_dro:joint_cc} to exceed its limit using $1000$ new samples. 
We computed these samples using a normal distribution truncated on $[\underline{\xi}_j, \overline{\xi}_j],\ j=1,...,D$ with zero mean and standard deviation $S^{\rm oos}_j = S_j + S_j^{\rm pert}(\epsilon_j)$, where $S_j$ is defined above and $S_j^{\rm pert}(\epsilon_j)$ is such that for a random variable $X\sim\mathcal{N}\big(0,S_j^{\rm pert}(\epsilon_j)\big)$ it holds that $\mathbb{E}\norm{X}_1 = \epsilon_j$ as per \cref{eq:wasserstein_additive}. 
We provide additional results for cases including $\epsilon_j=0$, i.e., when the decision maker assumes that the available data samples accurately reflect the true underlying distribution. 
We observe in Table~\ref{tab:oos_joint_cc} that all empirical violation probabilities for cases with $\epsilon_j > 0$ are smaller than the $5\%$ target. This is expected given that the CVaR produces a more conservative approximation of the initial chance constraint. For all $\big\{(\epsilon_1, \epsilon_2)|\epsilon_2\ge 0.1\big\}$ the empirical violation probability is zero because the model relied on the robust solution due to low data quality (i.e., $\lambda_1^{\rm cc}= \lambda_2^{\rm cc}=0$, see Table~II).
These results also corroborate our observation that dataset $j=2$ has an overall greater impact on the decision than $j=1$.
For some cases with $\epsilon_j = 0$ we observe instances where the empirical constraint violation probability exceeds the 5\% target. This indicates that the limited samples available in $\widehat{\set{X}}_j$ were not sufficient to accurately represent the true underlying distribution.
}

\begin{table}[]
\color{black}
    \centering
    \caption{Empirical out-of-sample violation probability of the joint chance constraint from $N=1000$ samples.}
    \label{tab:oos_joint_cc}
    \begin{tabular}{c|cccc|c}
     \diagbox{$\epsilon_1$}{$\epsilon_2$}   & \textbf{1.0} & \textbf{0.1} & \textbf{0.005} & \textbf{0.001} & \textbf{0.0}\\
        \hline
        \textbf{1.0} & 0.0\% & 0.0\% & 1.5\% & 3.6\% & 4.5\%\\
        \textbf{0.1} & 0.0\% & 0.0\% & 1.7\% & 3.0\% & 4.2\%\\
        \textbf{0.005} & 0.0\% & 0.0\% & 1.0\% & 3.1\% & 4.0\%\\
        \textbf{0.001} & 0.0\% & 0.0\% & 2.4\% & 3.2\% & 4.5\%\\
        \hline
        \textbf{0.0} & 3.6\% & 0.3\% & 2.8\% & 6.1\% & 7.4\%\\
    \end{tabular}
\end{table}

\subsection{Scalability}
\label{ssec:scalability}

All model runs for the results in Section~\ref{ssec:results_and_analysis} were completed in \unit[0.017]{s} on average. 
As outlined in Section~\ref{sec:multi-source_dro}, under conditions \ref{cond:separable_cost_and_support} \textit{or} \ref{cond:standardized_data}, which often hold in technical applications such as OPF, the size of the multi-source DRO scales linearly with the number of uncertain resources and data samples.
Hence, our multi-source DRO OPF problems scales identically to similar data-driven approaches, e.g., \cite{arrigo2022wasserstein,guo2018data}, and lends itself to the application of performance enhancing solution algorithms such as branch-and-cut or column-and-constraint-generation techniques.

\section{Conclusion}
\label{sec:conlcusion}

Motivated by the proliferation of power system data collection and analysis, we have proposed a multi-source distributionally-robust optimization approach that facilitates data-driven optimization from various datasets with heterogeneous quality. 
To this end, we have discussed a data quality metric based on Wasserstein distance and shown its usefulness for practical use cases related to differential privacy, noisy data measurements, and other algorithmic data obfuscation techniques. 
Using this data quality metric, we have derived a general stochastic optimization framework that leverages multi-source distributionally robust optimization (MS-DRO) and have proven and discussed its central properties. 
The proposed MS-DRO method enables data valuation for a multitude of decision-making problems. As a case study, we have focused on an exemplary stochastic OPF problem, which enabled us to highlight useful model properties, e.g., a direct approach to compute the marginal value of data quality and physical interpretations for both data quality and its marginal value. 
Our case study case illustrated the various formulation components and corroborated our analytical results.

There are two major avenues for future work. First, we plan to showcase the MS-DRO on a larger scale application, e.g., in a smart distribution system with data streams from various vendors. The resulting data valuation for different operation scenarios should provide a model to analyze investments in data infrastructure. 
Second, we plan to provide further analysis and potential mechanism design approaches for data remuneration. 
{\color{black} In particular, we will study approaches to ensure truthful reporting of data quality.}

\bibliographystyle{IEEEtran}
\bibliography{literature}

\begin{thebibliography}{10}
\providecommand{\url}[1]{#1}
\csname url@samestyle\endcsname
\providecommand{\newblock}{\relax}
\providecommand{\bibinfo}[2]{#2}
\providecommand{\BIBentrySTDinterwordspacing}{\spaceskip=0pt\relax}
\providecommand{\BIBentryALTinterwordstretchfactor}{4}
\providecommand{\BIBentryALTinterwordspacing}{\spaceskip=\fontdimen2\font plus
\BIBentryALTinterwordstretchfactor\fontdimen3\font minus
  \fontdimen4\font\relax}
\providecommand{\BIBforeignlanguage}[2]{{%
\expandafter\ifx\csname l@#1\endcsname\relax
\typeout{** WARNING: IEEEtran.bst: No hyphenation pattern has been}%
\typeout{** loaded for the language `#1'. Using the pattern for}%
\typeout{** the default language instead.}%
\else
\language=\csname l@#1\endcsname
\fi
#2}}
\providecommand{\BIBdecl}{\relax}
\BIBdecl

\bibitem{pjm2020pmus}
\BIBentryALTinterwordspacing
{PJM}. (2020) {Phasor Measurement Unit (PMU) Placement Plan in RTEP Planning
  Process}. [Online]. Available:
  \url{https://pjm.com/-/media/committees-groups/committees/pc/2020/20200707/20200707-item-05a-pmu-rtep.ashx}
\BIBentrySTDinterwordspacing

\bibitem{siemens2019gridedge}
\BIBentryALTinterwordspacing
Siemens, ``The grid edge revolution: Innovative drivers towards net-zero
  energy,'' Tech. Rep., 2019. [Online]. Available:
  \url{new.siemens.com/global/en/company/topic-areas/smart-infrastructure/grid-edge/white-paper-grid-edge-net-zero-energy-drivers.html}
\BIBentrySTDinterwordspacing

\bibitem{agarwal2019marketplace}
A.~Agarwal \emph{et~al.}, ``A marketplace for data: An algorithmic solution,''
  in \emph{Proc. ACM Conf. on Econ. and Comp.}, 2019.

\bibitem{acharya2022false}
S.~Acharya \emph{et~al.}, ``False data injection attacks on data markets for
  electric vehicle charging stations,'' \emph{Adv. Appl. Energy}, vol.~7, 2022.

\bibitem{giaconi2021smart}
G.~Giaconi \emph{et~al.}, ``Smart meter data privacy,'' in \emph{Advanced Data
  Analytics for Power Systems}.\hskip 1em plus 0.5em minus 0.4em\relax
  Cambridge University Press, 2021.

\bibitem{roald2023power}
L.~A. Roald \emph{et~al.}, ``Power systems optimization under uncertainty: A
  review of methods and applications,'' \emph{Electric Power Systems Research},
  vol. 214, p. 108725, 2023.

\bibitem{hassan2020data}
A.~Hassan \emph{et~al.}, ``Data-driven learning and load ensemble control,''
  \emph{Electric Power Systems Research}, vol. 189, p. 106780, 2020.

\bibitem{peng2022markovian}
G.~Peng \emph{et~al.}, ``Markovian decentralized ensemble control for demand
  response,'' \emph{IEEE Control Systems Letters}, vol.~6, pp. 3050--3055,
  2022.

\bibitem{dall2017chance}
E.~Dall’Anese \emph{et~al.}, ``Chance-constrained ac optimal power flow for
  distribution systems with renewables,'' \emph{IEEE Trans. Power Syst.},
  vol.~32, no.~5, pp. 3427--3438, 2017.

\bibitem{mieth2018data}
R.~Mieth and Y.~Dvorkin, ``Data-driven distributionally robust optimal power
  flow for distribution systems,'' \emph{IEEE Control Systems Letters}, vol.~2,
  no.~3, pp. 363--368, 2018.

\bibitem{morales2021learning}
J.~M. Morales \emph{et~al.}, ``Learning the price response of active
  distribution networks for {TSO-DSO} coordination,'' \emph{IEEE Trans. Power
  Syst.}, vol.~37, no.~4, pp. 2858--2868, 2021.

\bibitem{li2017distributed}
P.~Li \emph{et~al.}, ``A distributed online pricing strategy for demand
  response programs,'' \emph{IEEE Trans. Smart Grid}, vol.~10, no.~1, pp.
  350--360, 2017.

\bibitem{mieth2019online}
R.~Mieth and Y.~Dvorkin, ``Online learning for network constrained demand
  response pricing in distribution systems,'' \emph{IEEE Trans. Smart Grid},
  vol.~11, no.~3, pp. 2563--2575, 2019.

\bibitem{tucker2020constrained}
N.~Tucker \emph{et~al.}, ``Constrained {Thompson} sampling for real-time
  electricity pricing with grid reliability constraints,'' \emph{IEEE Trans.
  Smart Grid}, vol.~11, no.~6, pp. 4971--4983, 2020.

\bibitem{guo2018data}
Y.~Guo \emph{et~al.}, ``Data-based distributionally robust stochastic optimal
  power flow—{Part I}: Methodologies,'' \emph{IEEE Trans. Power Syst.},
  vol.~34, no.~2, pp. 1483--1492, 2018.

\bibitem{arrigo2022wasserstein}
A.~Arrigo \emph{et~al.}, ``Wasserstein distributionally robust
  chance-constrained optimization for energy and reserve dispatch: An exact and
  physically-bounded formulation,'' \emph{European Journal of Operational
  Research}, vol. 296, no.~1, pp. 304--322, 2022.

\bibitem{morales2023prescribing}
J.~M. Morales \emph{et~al.}, ``Prescribing net demand for two-stage electricity
  generation scheduling,'' \emph{Operations Research Perspectives}, vol.~10, p.
  100268, 2023.

\bibitem{crozier2022data}
C.~Crozier \emph{et~al.}, ``Data-driven contingency selection for fast security
  constrained optimal power flow,'' in \emph{Proc. 17th Int. Conference on
  Probabilistic Methods Applied to Power Systems (PMAPS)}.\hskip 1em plus 0.5em
  minus 0.4em\relax IEEE, 2022.

\bibitem{fernandez2021inverse}
R.~Fern{\'a}ndez-Blanco \emph{et~al.}, ``Inverse optimization with kernel
  regression: Application to the power forecasting and bidding of a fleet of
  electric vehicles,'' \emph{Computers \& Operations Research}, vol. 134, p.
  105405, 2021.

\bibitem{munoz2020feature}
M.~Munoz \emph{et~al.}, ``Feature-driven improvement of renewable energy
  forecasting and trading,'' \emph{IEEE Trans. Power Syst.}, vol.~35, no.~5,
  pp. 3753--3763, 2020.

\bibitem{ren2018datum}
X.~Ren \emph{et~al.}, ``Datum: Managing data purchasing and data placement in a
  geo-distributed data market,'' \emph{IEEE/ACM Trans. Networking}, vol.~26,
  no.~2, pp. 893--905, 2018.

\bibitem{bessa2018data}
R.~J. Bessa \emph{et~al.}, ``Data economy for prosumers in a smart grid
  ecosystem,'' in \emph{Proc. Ninth International Conference on Future Energy
  Systems}, 2018, pp. 622--630.

\bibitem{goncalves2020towards}
C.~Goncalves \emph{et~al.}, ``Towards data markets in renewable energy
  forecasting,'' \emph{IEEE Trans. Sust. Energy}, vol.~12, no.~1, pp. 533--542,
  2020.

\bibitem{han2021monetizing}
L.~Han \emph{et~al.}, ``Monetizing customer load data for an energy retailer: A
  cooperative game approach,'' in \emph{Proc. IEEE PowerTech}.\hskip 1em plus
  0.5em minus 0.4em\relax IEEE, 2021.

\bibitem{le2020ethical}
G.~Le~Ray and P.~Pinson, ``The ethical smart grid: Enabling a fruitful and
  long-lasting relationship between utilities and customers,'' \emph{Energy
  Policy}, vol. 140, p. 111258, 2020.

\bibitem{han2022trading}
L.~Han \emph{et~al.}, ``Trading data for wind power forecasting: A regression
  market with lasso regularization,'' \emph{Electric Power Systems Research},
  vol. 212, p. 108442, 2022.

\bibitem{pinson2022regression}
P.~Pinson \emph{et~al.}, ``Regression markets and application to energy
  forecasting,'' \emph{{TOP}}, vol.~30, no.~3, pp. 533--573, 2022.

\bibitem{raja2023market}
A.~A. Raja \emph{et~al.}, ``A market for trading forecasts: A wagering
  mechanism,'' \emph{International Journal of Forecasting}, 2023.

\bibitem{xie2022robust}
R.~Xie \emph{et~al.}, ``Robust scheduling with improved uncertainty sets via
  purchase of distributed predictive information,'' \emph{arXiv:2210.00291},
  2022.

\bibitem{vdvorkin2020differentially}
V.~Dvorkin \emph{et~al.}, ``Differentially private optimal power flow for
  distribution grids,'' \emph{IEEE Trans. Power Syst.}, vol.~36, no.~3, 2020.

\bibitem{brown2016characterizing}
M.~Brown \emph{et~al.}, ``Characterizing and quantifying noise in {PMU} data,''
  in \emph{Proc. IEEE Power and Energy Society General Meeting (PESGM)}.\hskip
  1em plus 0.5em minus 0.4em\relax IEEE, 2016, pp. 1--5.

\bibitem{esfahani2018data}
P.~M. Esfahani and D.~Kuhn, ``Data-driven distributionally robust optimization
  using the wasserstein metric: Performance guarantees and tractable
  reformulations,'' \emph{Math. Programming}, vol. 171, no. 1-2, 2018.

\bibitem{awasthi2022distributionally}
P.~Awasthi \emph{et~al.}, ``Distributionally robust data join,''
  \emph{arXiv:2202.05797}, 2022.

\bibitem{alvarez2012similarity}
P.~C. Alvarez-Esteban \emph{et~al.}, ``Similarity of samples and trimming,''
  2012.

\bibitem{kuhn2019wasserstein}
D.~Kuhn \emph{et~al.}, ``Wasserstein distributionally robust optimization:
  Theory and applications in machine learning,'' in \emph{Operations Research
  \& Management Science in the Age of Analytics}.\hskip 1em plus 0.5em minus
  0.4em\relax INFORMS, 2019, pp. 130--166.

\bibitem{duan2018distributionally}
C.~Duan \emph{et~al.}, ``Distributionally robust chance-constrained approximate
  ac-opf with wasserstein metric,'' \emph{IEEE Transactions on Power Systems},
  vol.~33, no.~5, pp. 4924--4936, 2018.

\bibitem{gopstein2021nist}
A.~Gopstein \emph{et~al.}, \emph{{NIST Framework and Roadmap for Smart Grid
  Interoperability Standards, Release 4.0}}.\hskip 1em plus 0.5em minus
  0.4em\relax Department of Commerce. National Institute of Standards and
  Technology, 2021.

\bibitem{molina2010private}
A.~Molina-Markham \emph{et~al.}, ``Private memoirs of a smart meter,'' in
  \emph{Proc. 2nd ACM workshop on embedded sensing systems for
  energy-efficiency in building}, 2010, pp. 61--66.

\bibitem{hassan2021privacy}
A.~Hassan \emph{et~al.}, ``Privacy-aware load ensemble control: A
  linearly-solvable mdp approach,'' \emph{IEEE Trans. Smart Grid}, vol.~13,
  no.~1, 2021.

\bibitem{dwork2014algorithmic}
C.~Dwork \emph{et~al.}, ``The algorithmic foundations of differential
  privacy,'' \emph{Foundations and Trends in Theoretical Computer Science},
  vol.~9, no. 3--4, pp. 211--407, 2014.

\bibitem{panaretos2019statistical}
V.~M. Panaretos and Y.~Zemel, ``Statistical aspects of wasserstein distances,''
  \emph{{Annual Review of Statistics and its Application}}, vol.~6, 2019.

\bibitem{wasserman2010statistical}
L.~Wasserman and S.~Zhou, ``A statistical framework for differential privacy,''
  \emph{Journal of the American Statistical Association}, vol. 105, no. 489,
  pp. 375--389, 2010.

\bibitem{kursawe2011privacy}
K.~Kursawe \emph{et~al.}, ``Privacy-friendly aggregation for the smart-grid,''
  in \emph{Proc. International Symposium on Privacy Enhancing Technologies
  Symposium}.\hskip 1em plus 0.5em minus 0.4em\relax Springer, 2011, pp.
  175--191.

\bibitem{mieth2020risktrading}
R.~Mieth \emph{et~al.}, ``Risk trading in a chance-constrained stochastic
  electricity market,'' \emph{IEEE Control Systems Letters}, vol.~5, no.~1, pp.
  199--204, 2020.

\bibitem{arrigo2022embedding}
A.~Arrigo \emph{et~al.}, ``Embedding dependencies between wind farms in dist.
  robust optimal power flow,'' \emph{IEEE Trans. Power Syst.}, 2022.

\bibitem{hassan2019stochastic}
A.~Hassan \emph{et~al.}, ``Stochastic and distributionally robust load ensemble
  control,'' \emph{IEEE Trans. Power Syst.}, vol.~35, no.~6, pp. 4678--4688,
  2020.

\bibitem{esteban2021distributionally}
A.~Esteban-P{\'e}rez and J.~M. Morales, ``Distributionally robust optimal power
  flow with contextual information,'' \emph{European Journal of Operational
  Research}, vol. 306, no.~3, pp. 1047--1058, 2023.

\bibitem{van2021data}
B.~P. Van~Parys \emph{et~al.}, ``From data to decisions: Distributionally
  robust optimization is optimal,'' \emph{Management Science}, vol.~67, no.~6,
  2021.

\bibitem{mieth2019risk}
R.~Mieth \emph{et~al.}, ``Risk-and variance-aware electricity pricing,''
  \emph{Electric Power Systems Research}, vol. 189, p. 106804, 2020.

\bibitem{roveto2020co}
M.~Roveto \emph{et~al.}, ``Co-optimization of var and cvar for data-driven
  stochastic demand response auction,'' \emph{IEEE Control Systems Letters},
  vol.~4, no.~4, pp. 940--945, 2020.

\bibitem{rockafellar2000optimization}
R.~T. Rockafellar \emph{et~al.}, ``Optimization of conditional value-at-risk,''
  \emph{{Journal of Risk}}, vol.~2, pp. 21--42, 2000.

\bibitem{kazempour2018stochastic}
J.~Kazempour \emph{et~al.}, ``A stochastic market design with revenue adequacy
  and cost recovery by scenario,'' \emph{IEEE Trans. Power Syst.}, vol.~33,
  no.~4, pp. 3531--3545, 2018.

\bibitem{fang2019introducing}
X.~Fang \emph{et~al.}, ``Introducing uncertainty components in locational
  marginal prices for pricing wind power and load uncertainties,'' \emph{IEEE
  Trans. Power Syst.}, vol.~34, no.~3, pp. 2013--2024, 2019.

\bibitem{mieth2019distribution}
R.~Mieth and Y.~Dvorkin, ``Distribution electricity pricing under
  uncertainty,'' \emph{IEEE Trans. Power Syst.}, 2019.

\bibitem{kuang2018pricing}
X.~Kuang \emph{et~al.}, ``Pricing chance constraints in electricity markets,''
  \emph{IEEE Trans. Power. Syst.}, vol.~33, no.~4, pp. 4634--4636, 2018.

\bibitem{matpowercase5}
\BIBentryALTinterwordspacing
{MATPOWER}. (2014) {CASE5 Power flow data}. [Online]. Available:
  \url{https://matpower.org/docs/ref/matpower5.0/case5.html}
\BIBentrySTDinterwordspacing

\bibitem{dvorkin2015uncertainty}
Y.~Dvorkin \emph{et~al.}, ``Uncertainty sets for wind power generation,''
  \emph{IEEE Trans. Power Syst.}, vol.~31, no.~4, pp. 3326--3327, 2015.

\bibitem{lubin2023jump}
M.~Lubin \emph{et~al.}, ``Jump 1.0: Recent improvements to a modeling language
  for mathematical optimization,'' \emph{Mathematical Programming Computation},
  2023, in press.

\bibitem{git_msdro_opf}
\BIBentryALTinterwordspacing
R.~Mieth. (2023) {Multi-source DRO OPF Code Supplement}. [Online]. Available:
  \url{https://github.com/mieth-robert/msdro_opf_public}
\BIBentrySTDinterwordspacing

\bibitem{ordoudis2021energy}
C.~Ordoudis \emph{et~al.}, ``Energy and reserve dispatch with distributionally
  robust joint chance constraints,'' \emph{Operations Research Letters},
  vol.~49, no.~3, pp. 291--299, 2021.

\end{thebibliography}

\appendix
\setcounter{equation}{0}
\renewcommand{\theequation}{\thesection.\arabic{equation}}

\subsection{Proof of Theorem~\ref{th:multi-source-wasserstein-dro}}
\label{ax:proof_of_th_wassertein_dro_general}

The multi-source Wasserstein DRO problem is given as
\begin{subequations}
\begin{alignat}{2}
   \inf_{\bm{x} \in X}  \sup_{\mathcal{Q}\in\mathcal{P}(\Xi)}& \int_\Xi c(\bm{x},\bm{\xi}_1,..., \bm{\xi}_D) d\mathcal{Q}(\bm{\xi}_1,..., \bm{\xi}_D) \hspace{-2cm}&&\\
    \rm{s.t.}\quad &
    P_{j\#}\mathcal{Q}=\mathcal{Q}_j,\quad&& j=1,...,D\\
    & W_{p}^{p}(\widehat{\mathcal{P}}_j, \mathcal{Q}_j) \leq \epsilon_j,\quad&&
    j=1,...,D \label{seq:marginal_balls}
\end{alignat}
\label{eq:ms_wasserstein_dro_complete}
\end{subequations}
Constraints \cref{seq:marginal_balls} can be replaced with 
\begin{equation}
    \int_{\Xi_j \times \Xi_j}\norm{\bm{\xi}' - \bm{\xi}_j}^p \pi_j(d\bm{\xi}', d\bm{\xi}_j) \leq \epsilon_j,\quad j=1,...,D,
    \label{eq:marginal_balls_projecton}
\end{equation}
where $\Xi_j:= {\rm proj}_{j}(\Xi)$, $\pi_{j}$ is a feasible coupling between $\widehat{\mathcal{P}}_j$ and $\mathcal{Q}_j$, and ${\rm proj}_{j}(\Xi)$ stands for the projection of the support set $\Xi$ onto the $j$-th coordinate axis.
Now, consider the product (empirical) measure $\widehat{\mathcal{P}}= \prod_{j=1}^D \widehat{\mathcal{P}}_j$, which allows to rewrite \cref{eq:marginal_balls_projecton} as
\begin{align}
   \int_{\Xi \times \Xi}\norm{\bm{\xi}'_{j} - \bm{\xi}_j}^p \Pi(d\bm{\xi}'_1,..., d\bm{\xi}'_D, d\bm{\xi}_1,..., d\bm{\xi}_D) \leq \epsilon_j, \nonumber \\ \hfill j=1,...,D,
\end{align}
where $\Pi$ is a feasible coupling between $\widehat{\mathcal{P}}$ and $\mathcal{Q}$.
We now can reformulate the inner supremum of \cref{eq:ms_wasserstein_dro_complete} as 
\begin{subequations}\label{DRO_MDP2}
\begin{alignat}{2}
   \sup_{\{\mathcal{Q}_{\iota}\in\mathcal{P}(\Xi),\ \iota \in \set{I}\}} 
   & \frac{1}{|\set{I}|}\!\sum_{\iota \in \set{I}}\!\int_\Xi\! c(\bm{\xi}_1, ..., \bm{\xi}_D\!) d\mathcal{Q}_{\iota}(\bm{\xi}_1,..., \bm{\xi}_D\!)\label{seq:split_objective1} \\
    \rm{s.t.}\ 
    & \frac{1}{|\set{I}|} \sum_{\iota \in \set{I}} \int_{\Xi}\norm{\bm{\xi}_{j} - \widehat{\bm{\xi}}_{\iota_j}}^p \mathcal{Q}_{\iota}(d\bm{\xi}_1, ..., d\bm{\xi}_D) \leq \epsilon_j, \nonumber \\
    & \hspace{3.5cm} j=1,...,D,\label{seq:marginal_balls2}
\end{alignat}%
\label{eq:general_multidata_wcexp}%
\end{subequations}%
where we dropped the dependency on $\bm{x}$ for leaner notation.
The remainder of the proof follows the proof of \cite[Theorem 4.2]{esfahani2018data}: (i) introducing dual variable $\lambda_j$ on \cref{seq:marginal_balls2}, (ii)~exchanging minimization and maximization operators, and (iii)~introducing epigraphical auxiliary variables $s_{\iota},\ \iota\in\set{I}$.

\subsection{Proof of Proposition~\ref{prop:mwsdro_sep_cost}}
\label{ax:proof_sep_cost}
Under Condition~\ref{cond:separable_cost_and_support} the expectation in \cref{seq:split_objective1} is 
\begin{equation}
\begin{aligned}
    &\frac{1}{|\set{I}|}\!\sum_{\iota \in \set{I}}\!\int_\Xi\! c(\bm{\xi}_1, ..., \bm{\xi}_D\!) d\mathcal{Q}_{\iota}(\bm{\xi}_1,..., \bm{\xi}_D\!)  \\
& = \sum_{j=1}^D\frac{1}{|\set{I}|}\!\sum_{\iota \in \set{I}}\!\int_{\Xi_j}\! c_j(\bm{\xi}_j) dP_{j\#}\mathcal{Q}_{\iota}(\bm{\xi}_j) \\
& = \sum_{j=1}^D\frac{1}{N_j}\!\sum_{i=1}^{N_j}\!\int_{\Xi_j}\! c_j(\bm{\xi}_j) d\mathcal{Q}_{ji}(\bm{\xi}_j)
\end{aligned}
\end{equation}
and, similarly, transportation budget constraints \cref{seq:marginal_balls2} are
\begin{equation}
    \frac{1}{N_j}\sum_{i=1}^{N_j}\int_{\Xi_j} \norm{\bm{\xi}_j - \widehat{\bm{\xi}}_{ji}}^{p}\mathcal{Q}_{ji}(d\bm{\xi}_j) \le \epsilon_j,\ j=1,...,D.
\end{equation}
The required result follows as in the final steps of Appendix~\ref{ax:proof_of_th_wassertein_dro_general}.

\subsection{Proof of Proposition~\ref{prop:mwsdro_stand_data}}
\label{ax:proof_stand_data}

Under Condition~\ref{cond:standardized_data} each sample $i'$ across all $j=1,...,D$ can be interpreted as a sample drawn from a $D$-dimensional distribution. Thus the expectation in \cref{seq:split_objective1} becomes:
\begin{equation}
    \frac{1}{N'}\sum_{i'=1}^{N'} \int_{\Xi} c(\bm{x}, \bm{\xi})\mathcal{Q}_{i'}(d\bm{\xi})
\end{equation}
and transportation budget constraints \cref{seq:marginal_balls2} become
\begin{equation}
    \frac{1}{N'} \sum_{i'=1}^{N'} \int_{\Xi} \norm{\bm{\xi}_j - \widehat{\bm{\xi}}_{ji'}}^p \mathcal{Q}_{i'}(d\bm{\xi}) \le \epsilon_j,\ j=1,...,D.
\end{equation}
The required result follows as in the final steps of Appendix~\ref{ax:proof_of_th_wassertein_dro_general}.

\subsection{Derivation details of objective reformulation }
\label{ax:wc_obj_derivation}
In \cref{eq:compact_dcopf_dro} we have the separable cost function 
\begin{equation}
    c(\bm{x}, \bm{\xi}) = \dprod{\bm{c}^{\rm A}}{-\bm{A}\bm{\xi}} = \sum_{j=1}^D \sum_{g=1}^G -c_g^{\rm A} a_{gj} \xi_{j}. 
\end{equation}
As a result, constraints in \cref{prop1:constraint} become
\begin{equation}
    s_{ij} \ge \sup_{\xi_j \in[\underline{\xi}_j, \overline{\xi}_j]} \sum_{g=1}^G -c_g^{\rm A} a_{gj} \xi_{j}  - \lambda_j|\xi_j - \widehat{\xi}_{ji}|,
\end{equation}
which can be reformulated into \cref{eq:opf_wcexp}.

\subsection{Proof of Proposition~\ref{prop:epsilon_analysis}}
\label{ax:proof_epsilon_analysis}
Consider the epigraph formulation in \cref{eq:opf_wcexp} and define $c'_j=-\sum_{g=1}^G c_g^{\rm A} \alpha_{gj}$ as a shorthand. 
As discussed in Section~\ref{ssec:data_quality_interpretaion}, the largest sum of $s^{co}_{ji}$ is attained if $\lambda_j^{\rm co}=0$, in which case $s^{co}_{ji} = c_j'\underline{\xi}_j,\ \forall i$ per \cref{opf_wcexp:epi_l}, i.e., the worst case, and the smallest sum of $s^{co}_{ji}$ is attained if $\lambda_j^{\rm co}=c_j'$, in which case $s^{co}_{ji} = c_j'\widehat{\xi}_j,\ \forall i$ as \cref{opf_wcexp:epi_u,opf_wcexp:epi_l,opf_wcexp:epi_av} are equal.
As \cref{opf_wcexp:epi_u} is always dominated by \cref{opf_wcexp:epi_l,opf_wcexp:epi_av}, write $s^{co}_{ji} = c_j'\underline{\xi}_j + mc_j'(\underline{\xi}_j - \widehat{\xi}_{ji})$ where we set $\lambda_j^{\rm co} = m c_j'$ with $m=[0,1]$. Using this definition the relevant part of \cref{opf_wcexp:objective} can be written as 
$mc_j' \epsilon_j + (N')^{-1}\sum_{i=1}^{N'}(c_j'\underline{\xi}_j - mc_j'(\widehat{\xi}_{ji} -\underline{\xi}_j)
= mc_j' \big[\epsilon_j - (N')^{-1}\sum\nolimits_{i=1}^{N'}(\widehat{\xi}_{ji} -\underline{\xi}_j)\big] + c_j'\underline{\xi}_j$.
If the term in rectangular parentheses is positive (case 1 in Proposition~\ref{prop:epsilon_analysis}), $m$ should be minimized leading to $m=0 \Rightarrow \lambda_j^{\rm co}=0$. 
If term (A) is negative (case 2 in Proposition~\ref{prop:epsilon_analysis}), $m$ should be maximized leading to $m=1 \Rightarrow \lambda_j^{\rm co}=c_j'$.

\end{document}